\documentclass{article}
\usepackage{amsmath, amssymb,amsfonts}
\usepackage{color}
\catcode`\@=11 \@addtoreset{equation}{section}

\catcode`\@=12

\newtheorem{Theorem}{Theorem}
\newtheorem{Proposition}{Proposition}[section]
\newtheorem{Lemma}{Lemma}[section]
\newtheorem{Corollary}{Corollary}[section]

\newcommand{\bTheorem}[1]{
\begin{Theorem} \label{T#1} }
\newcommand{\eT}{\end{Theorem}}

\newcommand{\bProposition}[1]{
\begin{Proposition} \label{P#1}}
\newcommand{\eP}{\end{Proposition}}

\newcommand{\bLemma}[1]{
\begin{Lemma} \label{L#1} }
\newcommand{\eL}{\end{Lemma}}

\newcommand{\bCorollary}[1]{
\begin{Corollary} \label{C#1} }
\newcommand{\eC}{\end{Corollary}}

\newcommand{\bFormula}[1]{
\begin{equation} \label{#1}}
\newcommand{\eF}{\end{equation}}

\newcommand{\Bl}{\left\{}
\newcommand{\Br}{\right\}}
\newcommand{\pa}{{\partial}}
\newcommand{\na}{{\nabla}}

\newcommand{\eps}{{\varepsilon}}

\newcommand{\va}[1]{\left<#1\right>}

\newcommand{\DC}{C^\infty_c}

\newcommand{\vr}{\varrho}

\newcommand{\vre}{\vr_\ep}

\newcommand{\vue}{\vu_\ep}

\newcommand{\vu}{\vc{u}}
\newcommand{\vc}[1]{{\bf #1}}
\newcommand{\qed}{\bigskip \rightline {Q.E.D.} \bigskip}
\newcommand{\Div}{{\rm div}_x}
\newcommand{\curl}{{\rm curl}_x}
\newcommand{\Divh}{{\rm div}_h}
\newcommand{\Grad}{\nabla_x}

\newcommand{\tn}[1]{\mbox {\F #1}}
\newcommand{\dx}{{\rm d} {x}}
\newcommand{\dt}{{\rm d} t }

\newcommand{\dxdt}{\dx \ \dt}
\newcommand{\intO}[1]{\int_{\Omega} #1 \ \dx}

\newcommand{\bProof}{{\bf Proof: }}

\newcommand{\ep}{\varepsilon}

\font\F=msbm10 scaled 1000
\newcommand{\R}{\mbox{\F R}}

\date{}
\begin{document}

%%%%%%%%%%%%%%%%%%%%%%%%%%%%%%%%

%%%%%%%%%%%%%%%%%%%%%%%%%%%%%%%%

\title{Multi-scale analysis of compressible viscous and rotating fluids}
\author{Eduard Feireisl\thanks{The work was supported by Grant 201/09/
0917 of GA \v CR as a part of the general research programme of the
Academy of Sciences of the Czech Republic, Institutional Research
Plan AV0Z10190503.} \and Isabelle Gallagher\thanks{The work was partially supported by
the A.N.R grant ANR-08-BLAN-0301-01 "Mathoc\'ean", as well as  the Institut Universitaire de France.}\and David Gerard-Varet\thanks{The work was partially supported by
the A.N.R grant ANR-08-JCJC-0104 "RUGO".}
\and Anton\' \i n Novotn\' y} \maketitle

\bigskip

\centerline{Institute of Mathematics of the Academy of Sciences of the Czech Republic}
\centerline{\v Zitn\' a 25, 115 67 Praha 1, Czech Republic}

\medskip

\centerline{Institut de Math\' ematiques UMR 7586, Universit\' e Paris Diderot}
\centerline{175, Rue de Chevaleret, 75013 Paris, France}

\medskip

\centerline{IMATH Universit\' e du Sud Toulon-Var}
\centerline{BP 20132, 83957 La Garde, France}

\begin{abstract}
We study a singular limit for the compressible Navier-Stokes system when the Mach  and Rossby numbers are proportional to certain powers of a small parameter $\ep$. If the Rossby number dominates the Mach number, the limit problem is represented by the 2-D incompressible Navier-Stokes system describing the horizontal motion of vertical averages of the velocity field.
If they are of the same order then the limit problem turns out to be a linear, 2-D equation with a unique radially symmetric solution.
The effect of the centrifugal force is taken into account.
\end{abstract}

\section{Introduction}

Rotating fluid systems appear in many applications of fluid mechanics, in particular in
models of atmospheric and geophysical flows, see the monograph \cite{CDGG}. Earth's rotation, together
with the influence of gravity and the fact that atmospheric Mach number is typically very small, give rise to a
large variety of singular limit problems, where some of these characteristic numbers become large or tend to zero,
see Klein
\cite{klein3}. We consider a simple situation, where the Rossby number is proportional to a small
parameter~$\ep$, while the Mach number behaves like~$\ep^m$,  {with $m \geq 1$.
Scaling of this type with various choices of~$m$ appears, for instance, in meteorological models
(cf. \cite[Section 1.3]{klein3}).}

\medskip

\noindent
Neglecting the influence of the temperature we arrive at the
following scaled \emph{Navier-Stokes system} describing the time evolution of the fluid density
$\vr = \vr(t,x)$ and the velocity field $\vu = \vu(t,x)$:
\bFormula{i1}
\partial_t \vr + \Div (\vr \vu) = 0,
\eF
\bFormula{i2}
\partial_t (\vr \vu) + \Div (\vr \vu \otimes \vu) + \frac{1}{\ep} (\vc{b} \times \vr \vu) +
\frac{1}{\ep^{2m}} \Grad p(\vr) = \Div \tn{S}(\Grad \vu) + \frac{1}{\ep^2} \vr \Grad G,
\eF
where $p$ is the pressure, and $\tn{S}$ is the viscous stress tensor determined by
Newton's rheological law
\bFormula{i2a}
\tn{S}(\Grad \vu) = \mu \Big( \Grad \vu + \Grad^t \vu - \frac{2}{3} \Div \vu \tn{I} \Big),\ \mu > 0.
\eF
{For the sake of simplicity, we have omitted possible influence of the so-called bulk viscosity component in the viscous stress. }

\medskip

\noindent
We consider a very simple geometry of the underlying physical space $\Omega \subset\R^3$,
namely $\Omega$ is an infinite slab,
\[
\Omega = \R^2 \times { (0,1)}.
\]
Moreover, to eliminate entirely the effect of the boundary on the motion, we prescribe the \emph{complete slip} boundary conditions for the velocity field:
\bFormula{i2b}
\vc{u} \cdot \vc{n}  |_{\partial \Omega} = 0,\
[\tn{S} \vc{n}] \times \vc{n} |_{\partial \Omega} = 0,
\eF
{where $\vc{n}$ denotes the outer normal vector to the boundary. Note that the more common
\emph{no-slip} boundary condition
\[
\vc{u}|_{\partial \Omega} = 0
\]
would yield a trivial result in the asymptotic limit, namely $\vc{u} \to 0$ for $\ep \to 0$. On the other hand,
the so-called \emph{Navier's boundary condition}
\bFormula{i2bNavier}
\vc{u} \cdot \vc{n}  |_{\partial \Omega} = 0,\ \beta \vc{u}_{\rm tan}  +
[\tn{S} \vc{n}]_{\rm tan} |_{\partial \Omega} = 0, \ \beta > 0,
\eF
gives rise to a friction term in the limit system known as Ekman's pumping, see Section \ref{m} below.
 }

\medskip

\noindent
As is well known (see Ebin \cite{EB}), the boundary conditions (\ref{i2b}) may be conveniently reformulated in terms of geometrical restrictions imposed on the state variables that are
\emph{periodic} with respect to the vertical variable $x_3$. More specifically, we take
\bFormula{i2c}
\Omega = \R^2 \times {\cal T}^1,
\eF
where ${\cal T}^1 = [-1,1]|_{ \{ -1,1 \} }$ is a one-dimensional torus, on which
the fluid density $\vr$ as well as the horizontal component of the velocity
$\vu_h = [u^1, u^2]$ { are extended to be even in $x_3$, while the vertical component $u^3$ is taken odd}:
\bFormula{sym}
\begin{array}{c}
\vr(x_1,x_2, - x_3) =  \vr(x_1,x_2,x_3),\
u^i(x_1,x_2, - x_3) = u^i(x_1,x_2, x_3), \ i=1,2,\\ \\
u^3(x_1,x_2, - x_3) =  - u^3(x_1,x_2, x_3).
\end{array}
\eF
Finally,
we assume that the rotation axis is parallel to $x_3$, namely
$\vc{b} = [0,0,1]$, and set $\Grad G \approx \Grad |x_h|^2$ - the associated centrifugal force, where we have written $x_h = [x^1,x^2]$.

\medskip

\noindent
As shown in \cite{CDGG}, \emph{incompressible} rotating fluids stabilize to a 2D motion
described by the vertical averages of the velocity provided the Rossby number $\ep$ is small enough. Besides, the stabilizing effect of rotation has been exploited by many authors, see e.g. Babin, Mahalov and Nicolaenko \cite{BaMaNi2}, \cite{BaMaNi1}. On the other hand,
compressible fluid flows in the low Mach number regime behave like the incompressible ones, see Klainerman and Majda \cite{KM1}, Lions and Masmoudi \cite{LIMA1}, among many others. Thus, at { at least for $m \gg1$}, solutions of the scaled system (\ref{i1}), (\ref{i2}) are first rapidly driven to incompressibility and then stabilize to a purely horizontal motion as $\ep \to 0$. On the other hand, the above mentioned scenario
changes completely if $m=1$. {In this case the speed of rotation and incompressibility
act on the same scale. Accordingly, the limit behaviour of the fluid is described by a single (linear) equation, see Section~\ref{m}. Note that such a picture is
in sharp contrast with \cite{FeGaNo}, where the effect of the centrifugal force is neglected. }

\medskip

\noindent
However, a rigorous justification of the above programme
is hampered by serious mathematical difficulties:

\begin{itemize}

\item The main issue in the low Mach number limit, at least in the case of the so-called
ill prepared initial data, is the presence of rapidly oscillating acoustic waves, cf.
Desjardins and Grenier \cite{DesGre}, Desjardins et al. \cite{DGLM}, Lions and Masmoudi
\cite{LIMA1}. Similarly to \cite{DesGre}, given the geometry of the spatial domain
$\Omega$, we may expect a local decay of the acoustic energy as a result of dispersive effects. Unfortunately, however, the fluid is driven by the centrifugal force that
becomes large for~$|x_h| \to \infty$. Specifically, we have
$G \approx \ep^{-2m}$ on the sphere of the radius $\ep^{-m}$, whereas the speed of sound in the fluid is proportional to $\ep^{-m}$. In other words, the centrifugal force changes effectively propagation of acoustic waves and this effect cannot be neglected, not even on compact subsets of the physical domain.

\item The dispersive estimates of Strichartz' type exploited in \cite{CDGG} cannot be
used in the present setting as the acoustic waves, represented by the gradient component of the velocity field, remain large for~$|x_h| \to \infty$.

\end{itemize}
Our approach is based on combination of dispersive estimates for acoustic waves with
the local method developed   in \cite{GallSR}. {As already mentioned, we
focus on two qualitatively different situations: one where $m \gg 1$, and one where~$m=1$.}

\begin{itemize}
\item {\textsc{ The case~$m \gg 1$.}}

In order to eliminate the effect of the centrifugal force, we compute the exact rate of the local decay of acoustic energy,
here proportional to $\ep^{m/2}$, by adapting the argument of Metcalfe \cite{Metca} ({cf. also
D'Ancona and Racke \cite{DancRac}}, Smith and Sogge \cite{SmSo}). Accordingly, the associated acoustic equation can be localized to balls of radii $\ep^{-\alpha}$ {for a certain $\alpha > m/2$, on which} $G(x) \approx
\ep^{-2 \alpha}$. Having established local decay of acoustic waves, we use the method developed in \cite{GallSR}, based on cancellation of several quantities in the convective term, similar to the local approach of Lions and Masmoudi \cite{LIMA6}.

\item {\textsc{ The case~$m = 1 $.}}

In this case both the high rotation and weak compressibility limits occur at the same scale. We can compute immediately the limiting diagnostic equations, and, similarly to the previous case, careful analysis of  cancellations in the convective term allows to identifiy the limit system. {In contrast with the situation studied in \cite{FeGaNo},
the limit system is linear as a consequence of strong stratification caused by the centrifugal force.}

\end{itemize}
The paper is organized in the following way.
In Section \ref{p}, we
introduce the concept of finite energy weak solutions to the compressible Navier-Stokes system (\ref{i1} - \ref{i2b}) and recall their basic properties. In particular, we present uniform bounds on solutions independent of the scaling parameter $\ep \to 0$. The  main results concerning the singular limit of solutions to
(\ref{i1} - \ref{i2b}) are stated in Section \ref{p}. The anisotropic situation, when~$m \gg 1,$  is {analyzed in Section~\ref{proof1} by means of several steps: The easy part concerns the formal identification
of the limit system, while, in Section \ref{a}, the propagation of acoustic waves is studied as well as their local decay in $\Omega$. Finally, the limit systems for $m\gg 1$ is justified in Section \ref{c} by means of a careful analysis of the convective term. Finally, the isotropic case, when~$m = 1 $, is examined in Section~\ref{proof2}, where the proofs borrow several ingredients introduced in \cite{GallSR}.}

\section{Preliminaries and statement of the main results}
\label{p}

{
In this section we introduce the main hypotheses, recall some known results concerning existence
of solutions to the primitive system as well as the nowadays standard uniform bounds independent of the scaling
parameter, and, finally, formulate our main results.}

\subsection{Main assumptions}
Consider a family of solutions $\vre$, $\vue$ of the Navier-Stokes system (\ref{i1} - \ref{i2b}) in $(0,T) \times \Omega$, emanating from the initial data
\[
\vre(0, \cdot) = \vr_{0,\ep}, \ \vue(0, \cdot) = \vu_{0, \ep}.
\]
We assume that the initial data are \emph{ill-prepared}, specifically,
\bFormula{u1+}
\vr_{0,\ep} = \tilde \vr_\ep + \ep^m r_{0,\ep},
\eF
where $\tilde \vre$ is a solution of the associated \emph{static problem}:
\[
\Grad p(\tilde \vr_{\ep}) = \ep^{2(m-1)} \tilde \vr_\ep \Grad G \ \mbox{in}\ \Omega.
\]
Consequently,
\[
P(\tilde \vr_\ep) = \ep^{2(m-1)} G + {\rm const}, \ \mbox{where}\ P(\vr) = \int_1^\vr \frac{p'(z)}{z} \ {\rm d}z.
\]
{Furthermore,
we suppose that
\bFormula{HYP1}
p \in C^1[0,\infty) \cap C^2(0,\infty) , \ p'(\vr) > 0 \ \mbox{for all}\ \vr > 0, \
\lim_{\vr \to \infty} \frac{p'(\vr)}{\vr^{\gamma - 1}} = c > 0
\eF
for a certain $\gamma > 1$ specified below,} and that
\bFormula{u1--}
\begin{array}{c}
G \in W^{1,\infty}(\Omega),\ G(x) \geq  0, \ G(x_1,x_2,-x_3) = G(x_1,x_2,x_3),
\\ \\
 |\Grad G(x)| \leq c(1 + |x_h|)  \ \mbox{for all}\ x \in \Omega.
\end{array}
\eF
Finally, we normalize
\bFormula{stac}
P(\tilde \vr_\ep) = \ep^{2(m-1)} G,
\eF
noticing { that
\bFormula{stac2}
  \ \tilde \vre(x) \geq 1,\ \tilde \vr_\ep(x) \to 1
\ \ \mbox{for any}\ x \in \Omega
\ \mbox{as}\ \ep \to  0 \ \mbox{for} \ m > 1,
\eF
whereas
\bFormula{stac1}
\tilde \vre \equiv  \tilde \vr \ \mbox{is independent of} \ \ep \ \mbox{provided}
\ m = 1.
\eF
}

\subsection{Energy inequality}

Introducing the {\emph{relative entropy}}
\[
E(\vr, \tilde \vr) := H(\vr) - H'(\tilde \vr)(\vr - \tilde \vr) - H(\tilde \vr),\quad  H(\vr) := \vr \int_1^\vr \frac{p(z)}{z^2} {\rm d}z,
\]
we note that $H'(\vr) = P(\vr) + {\rm const}$; therefore we may assume that
\begin{equation}\label{u1}
\begin{aligned}
  \intO{
  \left( \frac{1}{2} \vre |\vue|^2 + \frac{1}{\ep^{2m}} E(\vre, \tilde \vr_\ep) \right)(\tau, \cdot)}
+ \int_0^\tau \intO{ \tn{S} (\Grad \vue) : \Grad \vue } \:{\rm d} \tau'&\\
  \quad  \leq \intO{\left( \frac{1}{2} \vr_{0,\ep} |\vu_{0,\ep}|^2 + \frac{1}{\ep^{2m}} E(\vr_{0,\ep}, \tilde \vr_\ep) \right)} &
\end{aligned}
\end{equation}
including implicitly the \emph{mass compatibility condition}
\[
\intO{ (\vre - \tilde \vre) } = 0.
\]
Now, in order to establish uniform bounds independent of the scaling parameter $\ep$,
the initial data specified in (\ref{u1+}) must be chosen in such a way
that the expression on the right-hand side of (\ref{u1}) remains bounded uniformly for $\ep \to 0$. Thus, if $\gamma \leq 2$ in (\ref{HYP1}), it is enough to assume that
\[
\{ r_{0,\ep} \}_{\ep > 0} \ \mbox{bounded in}\ L^2 \cap L^\infty (\Omega),\
\{ \sqrt{\tilde \vre} \vu_{0,\ep} \}_{\ep > 0} \ \mbox{bounded in}\ L^2 (\Omega; R^3).
\]
In general, we suppose that
\bFormula{u1++}
\begin{array}{c}
\left\{ {\tilde \vre}^{\frac{\gamma - 2}{2}} r_{0, \ep} \right\}_{\ep > 0}
\ \mbox{bounded in} \ L^2(\Omega),\
\{ r_{0,\ep} \}_{\ep > 0} \ \mbox{bounded in}\ L^2 \cap L^\infty (\Omega),\\ \\
\{ \sqrt{\tilde \vre} \vu_{0,\ep} \}_{\ep > 0} \ \mbox{bounded in}\ L^2 (\Omega; R^3).
\end{array}
\eF

\subsection{Finite energy weak solutions}
\label{few}

We say that $\vr$, $\vu$ is a finite energy weak solution of the Navier-Stokes system
(\ref{i1} - \ref{i2a}) in $(0,T) \times \Omega$, supplemented with the initial data
(\ref{u1+}), if:

\begin{itemize}
\item
the energy inequality
$$
\begin{aligned}
\intO{ \left( \frac{1}{2} \vr |\vu|^2 + \frac{1}{\ep^{2m}} E(\vr, \tilde \vr_\ep) \right)(\tau, \cdot) }
+ \int_0^\tau \intO{ \tn{S} (\Grad \vu) : \Grad \vu } \:{\rm d} \tau' 
& \\
\quad \leq\intO{ \left( \frac{1}{2} \vr_{0,\ep} |\vu_{0,\ep}|^2 + \frac{1}{\ep^{2m}} E(\vr_{0,\ep}, \tilde \vr_\ep) \right) }
\end{aligned}
$$
holds for a.a. $\tau \in (0,T)$;

\item equation (\ref{i1}) is satisfied in the sense of distributions, specifically,
\[
\int_0^T \intO{ \Big( \vr \partial_t \varphi + \vr \vu \cdot \Grad \varphi \Big) }\ \dt =
- \intO{ \vr_{0, \ep} \varphi(0, \cdot) }
\]
for any $\varphi \in \DC([0,T) \times \Omega)$;

\item the pressure $p$ is locally integrable in $[0,T) \times \Omega$, equation (\ref{i2}) holds in the sense of distributions:
\begin{equation}\label{i2weak}
\begin{aligned}
\int_0^T \intO{ \Big( \vr \vu \cdot \partial_t \varphi + \vr \vu \otimes \vu : \Grad
\varphi - \frac{1}{\ep} \vc{b} \times (\vr \vu) \cdot \varphi + \frac{1}{\ep^{2m}}
p(\vr) \Div \varphi \Big) } \ \dt & \\
= \int_0^T \intO{ \Big( \tn{S}(\Grad \vu) : \Grad \varphi - \frac{1}{\ep^2} \vr \Grad G \cdot \varphi \Big)}
\ \dt - \intO{ \vr_{0,\ep} \vu_{0,\ep} \cdot \varphi (0, \cdot) }&
\end{aligned}
\end{equation}
for any $\varphi \in \DC([0,T) \times \Omega; R^3)$.
\end{itemize}
Under hypothesis (\ref{HYP1}), the \emph{existence} of finite energy weak solutions can be established by the method developed by Lions \cite{LI4}, with the necessary modifications specified in \cite{FNP} in order to handle the physically relevant range of adiabatic exponents $\gamma > 3/2$.

\subsection{Uniform bounds}

\label{u}

In order to study the asymptotic behavior of solutions, we first establish \emph{uniform bounds} independent of the scaling parameter $\ep \to 0$. As a matter of fact, all of them follow from the energy inequality (\ref{u1}).
{Similarly to \cite[Chapter 5]{FEINOV}, }
{we introduce the \emph{essential} and \emph{residual} component of
a function $h$ as}
\[
h = h_{\rm ess} + h_{\rm res},
\]
where
\[
h_{\rm ess}(t,x) = h(t,x)
\ \mbox{for} \ (t,x) \ \mbox{such that}\ \vre(t,x) \in (1/2,2),\
h_{\rm ess}(t,x) = 0 \ \mbox{otherwise}.
\]
Now, {by virtue of (\ref{stac}),
\bFormula{odhad}
1 \leq \tilde \vre(x) \leq 1 + c(r) \ep^{{2(m-1-\alpha)} }
\ \mbox{for all} \ x \in B_{r/\ep^\alpha},\ 0 \leq \alpha \leq m-1,
\eF
where we have denoted
\[
B_R = \{ x \in \Omega \ |\ |x_h| \leq R \}.
\]
}
It follows directly from energy inequality (\ref{u1}) that when~$m>1+\alpha$
\begin{eqnarray}
\label{u1+++} {\rm ess} \sup_{t \in (0,T)} \left\|  \left[ \frac{\vre - \tilde \vr_{\ep} }{\ep^m} \right]_{\rm ess} \right\|_{L^2(B_{r/ \ep^\alpha})} \leq c(r) & \\ 
\label{u1++++}{\rm ess} \sup_{t \in (0,T)} \int_{B_{r/\ep^\alpha}} [\vre]^\gamma_{\rm res} \ \dx \leq \ep^{2m} c(r)& \\ 
\label {u14+}
{\rm ess} \sup_{t \in (0,T)} \int_{B_{r/\ep^\alpha}}  1_{\rm res} \ \dx \leq \ep^{2m} c
(r)& \\ 
\label{u14++}
{\rm ess} \sup_{t \in (0,T)}
\left| \{ x \in \Omega \ | \ \vre(t,x) \leq 1/2 \} \right| \leq c {\ep^{2m}}
\leq c& \\ 
\label{u15+}
{\rm ess} \sup_{t \in (0,T)} \| \sqrt{\vre} \vue \|_{L^2(\Omega;R^3)} \leq c&
\end{eqnarray}
and
\bFormula{u15++}
\int_0^T \int_\Omega \tn{S} (\Grad \vue) : \Grad \vue \ \dxdt \leq c.
\eF
In the case when~$m = 1$, the bounds~(\ref{u1+++})-(\ref{u14+}) should be replaced by
\begin{eqnarray}
\label{u1+++1} {\rm ess} \sup_{t \in (0,T)} \left\|  \left[ \frac{\vre - \tilde \vr_{\ep} }{\ep^m} \right]_{\rm ess} \right\|_{L^2(B_{r})} \leq c(r) & \\ 
\label{u1++++1}{\rm ess} \sup_{t \in (0,T)} \int_{B_{r}} [\vre]^\gamma_{\rm res} \ \dx \leq \ep^{2m} c(r)& \\ 
\label {u14+1}
{\rm ess} \sup_{t \in (0,T)} \int_{B_{r}}  1_{\rm res} \ \dx \leq \ep^{2m} c
(r).&
\end{eqnarray}
Finally, combining (\ref{u14++} - \ref{u15++}) with
a variant of Korn's inequality (see \cite[Theorem 10.17]{FEINOV}), we obtain
\bFormula{u16+}
\int_0^T \| \vue \|^2_{W^{1,2}(\Omega; R^3)} \ \dt  \leq c.
\eF
{All generic constants in the previous estimates are independent of the scaling parameter $\ep$.}

\subsection{Main results}
\label{m}

In what follows,
the symbol $\vc{H}$ denotes the standard Helmholtz projection onto the space of solenoidal functions in $\Omega$, {specifically,
\[
\widehat {\vc{H} [\vc{v}]} (\xi, k) = {\vc{v}} (\xi,k) - \frac{1}
{|\xi|^2 + k^2 }\Big[ \xi \Big( \xi \cdot \widehat {\vc{v}}_h (\xi,k)  + k \widehat{v}_3 (\xi,k) \Big) , k  \Big( \xi \cdot \widehat {\vc{v}}_h (\xi,k)  + k \widehat{v}_3 (\xi,k) \Big)\Big],
\]
\[
\vc{H}^\perp [\vc{v}] = \vc{v} - \vc{H}[\vc{v}],
\]
where the symbol $\widehat{\vc{v}}(\xi,k)$, $\xi \in R^2, k \in Z$, denotes the Fourier transform of
$\vc{v} = \vc{v}(x_h,x_3)$}. Similarly, the Laplace operator $\Delta$ is identified
through
\[
\widehat{ \Delta v } \approx  -(|\xi|^2 + k^2) \ \widehat v (\xi,k).
\]
Finally, we introduce the vertical average of a function $v$ as
\begin{equation}\label{defvertical}
\va{v}(x_h) = \frac{1}{|{\cal T}^1|} \int_{{\cal T}^1} v(x_h, x_3) \ {\rm d}x_3.
\end{equation}

{ \subsubsection{Multiscale limit ($m \gg 1$)} }

\bTheorem{m1}
Let the pressure $p$ and the potential of the driving force $G$ satisfy hypotheses
(\ref{HYP1}), (\ref{u1--}), with $\gamma > 3/2$. Let $\vre$, $\vue$ be a
finite energy weak solution of the Navier-Stokes system in $(0,T) \times \Omega$
belonging to the symmetry class (\ref{sym}), emanating from the initial data
(\ref{u1+}), (\ref{u1++}). In addition, suppose that
\[
m> 10
\]
and that
\[
\vu_{0,\ep} \to \vc{U}_0 \ \mbox{weakly in}\ L^2(\Omega;R^3).
\]
Then
\[
{\rm ess} \sup_{t \in (0,T)} \| \vre - 1 \|_{(L^2 + L^\gamma )(K)} \leq \ep^m c(K)
 \ \mbox{for any compact}\ K \subset \Omega,
\]
\[
\vue \to \vc{U} \ \mbox{weakly in} \ L^2(0,T; W^{1,2}(\Omega;R^3)),
\]
where $\vc{U} = [\vc{U}_h(x_h),0]$ is the unique solution to the 2D incompressible Navier-Stokes system
\bFormula{NS1}
{\rm div}_h \vc{U}_h = 0,
\eF
\bFormula{NS2}
\partial_t \vc{U}_h + {\rm div}_h (\vc{U}_h \otimes \vc{U}_h) + \nabla_h \Pi = \mu \Delta_h \vc{U}_h,
\eF
with the initial data
\[
\vc{U}_h(0,\cdot) = \Big[ \vc{H} \Big[  [\va{\vc{U}_0}_h , 0]  \Big] \Big]_h
\]
\eT

\medskip

\noindent
{\bf Remark \ref{m}.1}
\noindent
{\it
{A short inspection of the proof of Theorem \ref{Tm1} given in Section \ref{proof1} below reveals
that replacing the complete slip condition (\ref{i2b}) by the Navier's slip condition (\ref{i2bNavier})
would produce an extra term $\beta \vc{U}_h$ on the left-hand side of (\ref{NS2}) known as  Ekman's pumping.
}
}

\subsubsection{Stratified limit ($m=1$)}

\bTheorem{m2}
Let the pressure $p$ satisfy hypotheses
(\ref{HYP1}), with $\gamma > 3$ and let $G(x_h) = |x_h|^2$. Let $\vre$, $\vue$ be a
finite energy weak solution of the Navier-Stokes system in $(0,T) \times \Omega$
belonging to the symmetry class~(\ref{sym}), emanating from the initial data
(\ref{u1+}), (\ref{u1++}), where
\[
m = 1.
\]
In addition, suppose that  {
\[
r_{0,\ep} \to r_0 \ \mbox{weakly in}\ L^2(\Omega), \
\vu_{0,\ep} \to \vc{U}_0 \ \mbox{weakly in}\ L^2(\Omega;R^3).
\]
Then
\[
r_\ep \equiv \frac{\vre - \tilde \vr }{\ep} \to r \ \mbox{weakly-(*) in} \ L^\infty(0,T;L^2(K)) \ \mbox{for any compact}\ K \subset {\Omega},
\]
\[
\vue \to \vc{U} \ \mbox{weakly in}\ L^2(0,T; W^{1,2}(\Omega;R^3)),
\]
with
\begin{equation}\label{kernel1}
r = r(t, x_h) \:  \mbox{ radially symmetric }, \quad \vc{U} = [\vc{U}_h(t, x_h),0],
\end{equation}
and
\begin{equation}\label{kernel2}    \na_h (P'(\tilde \vr) r )  +  \vc{U}_h^\perp = 0.
\end{equation}
Moreover, the function $r$ satisfies
\bFormula{limiteq}
\partial_t \Big( r - {\rm div}_h (\tilde \vr  \na_h (P'(\tilde \vr) r)) \Big) + \Delta^2_h  \Big( P'(\tilde \vr) r \Big)  = 0.
\eF
 In addition,  the initial value $r(0)$ is the unique radially symmetric function satisfying the integral identity
\bFormula{INIT}
\int_{R^2} \Big( \tilde \vr \na_h \left( P'(\tilde \vr) r(0) \right) \cdot \na_h \psi \: + \: r(0) \, \psi \Big) {\rm d}x_h \:  = \:  \int_{R^2} \Big(   \left<  \tilde \vr \vc{U}_{0,h} \right>  \cdot \na_h^\perp \psi     + \left< r_0 \right> \psi \Big) {\rm d}x_h
\eF
}
for all radially symmetric $\psi = \psi(x_h) \in \DC(R^2)$.
\eT
\medskip
{\bf Remark \ref{m}.2}
\noindent
{\it
The initial condition \eqref{INIT} can be interpreted in polar coordinates as
$$
r(0, s ) - \frac{1}{s} \partial_s
\Big( s \tilde \vr \partial_s \Big( P'(\tilde \vr) r(0,s) \Big) \Big)
= \frac{1}{2 \pi s} \int_{ |x_h| = s } \Big(
{\rm curl}_h \left< \tilde \vr \vc{U}_{0,h} \right> + \left< r_0 \right>  \Big) {\rm dS}_{x_h}.
$$
where $s = | x_h|$,  provided all quantities are sufficiently smooth.}

\medskip

\noindent
{\bf Remark \ref{m}.3}
\noindent
{\it 
Note that (\ref{kernel2}) implies that $\nabla_h \tilde \vr \cdot \vc{U}_h = 0$, or, equivalently,
\[
\vc{U}_h \cdot x_h = 0, \ \mbox{meaning the limit velocity} \ \vc{U}_h \ \mbox{is tangent to the level sets} \ \big\{ |x_h| = {\rm const}
\big\},
\]
Moreover it follows from (\ref{kernel1}-\ref{kernel2}) that $|\vc{U}_h|$ is constant on~$\big\{|x_h| = {\rm const}\big\}$, and~${\rm div}_h \vc{U}_h = 0$.
}

\medskip

\noindent It is remarkable that the limit  equation (\ref{limiteq}) is linear for $m=1$, in sharp contrast with the homogeneous case $\tilde \vr = {\rm const}$ treated in \cite{FeGaNo}. This is related to the fact that the limit density $\tilde \vr$ is
stratified (non-constant). More precisely, the absence of nonlinearity is related to the smallness of the kernel of the penalized operator, defined by \eqref{kernel2}. This phenomenon had already been identified for fluids with variable rotation axis, see \cite{GallSR}. Much of the analysis that we shall follow to prove Theorem~\ref{Tm2} borrows to this last reference.

\medskip

\noindent
The rest of the paper is devoted to the proof of Theorems \ref{Tm1} and~\ref{Tm2}, in Sections~\ref{proof1} and~\ref{proof2} respectively.

\section{Anisotropic scaling: Proof of Theorem~\ref{Tm1}}

\label{proof1}

\subsection{Preliminary remarks}

We start with some simple observations that follow directly from the uniform bounds (\ref{u1+++} - \ref{u16+}). Clearly, relations (\ref{odhad}), (\ref{u1+++}), and (\ref{u1++++}) imply that
\[
{\rm ess} \sup_{t \in (0,T)} \| \vre - 1 \|_{(L^2 + L^\gamma )(K)} \leq \ep^m c(K)
 \ \mbox{for any compact}\ K \subset \Omega,
\]
while (\ref{u16+}) yields immediately
\[
\vue \to \vc{U} \ \mbox{weakly in}\ L^2(0,T; W^{1,2}(\Omega;R^3)),
\ \mbox{at least for a suitable subsequence.}
\]
Consequently, letting $\ep \to 0$ in (\ref{i1}) yields
\bFormula{M1}
\Div \vc{U} = 0 \ \mbox{a.a. in}\ (0,T) \times \Omega.
\eF
Moreover, it follows from equation (\ref{i2}) that
\[
\vc{H}[\vc{b} \times \vc{U}] = 0, \ \mbox{meaning,}\
\vc{b} \times \vc{U} = \Grad \Phi \ \mbox{for a certain potential} \ \Phi.
\]
Consequently, $\Phi$ and $\vc{U}_h$ are independent of the vertical coordinate $x_3$, and moreover, ${\rm div}_h \vc{U}_h = 0$. Since~$\vc{U}$ is solenoidal, we have
$\partial_{x_3} U^3 = 0$; whence, as $U^3$ has zero vertical mean,
\[
U^3 = 0, \ \vc{U} = [\vc{U}_h(t,x_h),0].
\]
Finally, in view of the uniform bounds established in (\ref{u1+++}), (\ref{u16+}), it is easy to pass to the limit in all terms appearing in the weak formulation of momentum equation (\ref{i2weak}), tested on $\varphi = [\phi (t,x_h),0]$,   with~${\rm div}_h \phi = 0$, with the exception of the convective term $\Div (\vre \vue \otimes \vue)$ that will be analyzed in the next two sections.

\subsection{Acoustic waves}
\label{a}

Since $m > 10$, there exists $\alpha$ such that
\bFormula{akuhyp}
1 + \frac{m}{2} < \alpha < \frac{3}{4} (m - 2)\cdotp
\eF
This choice of the parameter $\alpha$ will become clear in Section \ref{c}.
Moreover, we introduce a family of cut-off functions $\chi_\ep$ such that
\bFormula{cutoff}
\begin{array}{c}
\displaystyle
\chi_\ep \in \DC(R^2), \ 0 \leq \chi_\ep \leq 1,\ \chi_\ep (x_h) = 1 \ \mbox{for} \ |x_h| \leq
\frac{1}{\ep^\alpha}  \ , \\ \\ \chi_\ep (x_h) = 0 \ \mbox{for} \ |x_h| \geq \frac{2}{\ep^\alpha},\
|\Grad \chi_\ep(x_h) | \leq 2 \ep^\alpha \ \mbox{for}\ x_h \in R^2.
\end{array}
\eF

\medskip
\noindent
As the density becomes constant in the asymptotic limit,
the basic idea  is to
``replace'' $\vue \approx \chi_\ep \vre \vue$ and write
\[
\chi_\ep \vre \vue = \vc{H}[\chi_\ep \vre \vue] + \vc{H}^\perp [\chi_\ep \vre \vue] =
\va{\vc{H}[\chi_\ep \vre \vue]} + \Big( \vc{H}[\chi_\ep \vre \vue] - \va{\vc{H}[\chi_\ep \vre \vue]}\Big) + \Grad \Psi_\ep,
\]
where $\Grad \Psi_\ep = \vc{H}^\perp [\chi_\ep \vre \vue]$.
 
\medskip
\noindent
In the remaining part of this section,
we examine the asymptotic behavior of the \emph{acoustic potential}~$\Psi_\ep$.
More specifically, we show that $\Grad \Psi_\ep$ tends to zero on compact subsets of $\Omega$ and therefore becomes negligible in the asymptotic limit $\ep \to 0$.

\subsubsection{Acoustic equation}

Following Lighthill \cite{LIGHTHILL1}, \cite{LIGHTHILL2}, we rewrite the Navier-Stokes system
(\ref{i1}), (\ref{i2}) in the form:
$$
\begin{aligned}
&\ep^m \partial_t \left( \frac{\vre - \tilde \vre}{\ep^m} \right) + \Div(\vre \vue) = 0, \\
&\ep^m \partial_t (\vre \vue)  + p'(1) \Grad \left( \frac{\vre - \tilde \vre}{\ep^m} \right) =  -  \frac{1}{\ep^m} \Grad \Big( p(\vre) - p(\tilde \vr_{\ep}) - p'(1)
(\vre - \tilde \vr_{\ep}) \Big)\\
&\qquad +  \ep^m \Div \Big( \tn{S} (\Grad \vue) - (\vre \vue \times \vue)  \Big) - \ep^{m-1} \left( \vc{b} \times \vre \vue \right)
+ \ep^{2(m-1)}  \frac{ \vre - \tilde \vr_{\ep} }{\ep^m}  \Grad G  
\end{aligned}
$$
where, exactly as in Section \ref{few}, equations are understood in the sense of distributions in $(0,T) \times \Omega$. Furthermore, introducing new variables
\[
S_\ep = \chi_\ep \frac{\vre - \tilde \vre}{\ep^m}, \ \vc{m}_\ep =
\chi_\ep \vre \vue,
\]
where $\chi_\ep$ is the cut-off function specified in (\ref{cutoff}),
we arrive at the equation
\bFormula{aku1r}
\ep^m \partial_t S_\ep  + \Div \vc{m}_\ep = \Grad \chi_\ep \cdot (\vre \vue),
\eF
while~$\vc{m}_\ep $ satisfies
\bFormula{aku2r}
\begin{aligned}
&\ep^m \partial_t \vc{m}_\ep  + p'(1) \Grad S_\ep= - \ep^{m-1} \left( \vc{b} \times \vc{m}_\ep  \right)
+ \ep^{2(m-1)} \chi_\ep \frac{ \vre - \tilde \vr_{\ep} }{\ep^m}  \Grad G +
 p'(1) \Grad \chi_\ep \frac{\vre - \tilde \vre}{\ep^m} \\
& {}\quad +
\ep^m \Div \Big( \chi_\ep \tn{S} (\Grad \vue) - \chi_\ep (\vre \vue \times \vue) \Big) - \ep^m \Big( \tn{S} (\Grad \vue) - \vre \vue \times \vue \Big) \cdot \Grad \chi_\ep
\\
&\quad- \frac{1}{\ep^m} \Grad \left[ \chi_\ep \Big( p(\vre) - p(\tilde \vr_{\ep}) - p'(1)
(\vre - \tilde \vr_{\ep}) \Big) \right] + \frac{1}{\ep^m} \Big( p(\vre) - p(\tilde \vr_{\ep}) - p'(1)
(\vre - \tilde \vr_{\ep}) \Big) \Grad \chi_\ep.
\end{aligned}
\eF
\subsubsection{Uniform bounds}

Our next goal is to deduce uniform bounds on all quantities appearing in the
acoustic equation (\ref{aku1r}), (\ref{aku2r}). To begin, it
follows from (\ref{u1+++} - \ref{u14+}) that
\bFormula{U1}
\{ S_\ep \}_{\ep > 0} \ \mbox{is bounded in} \ L^\infty(0,T;L^2 + L^1 (\Omega)).
\eF
As a matter of fact, we have
\[
\{ S_\ep \}_{\ep > 0} \ \mbox{is bounded in} \ L^\infty(0,T;(L^2 + L^1 \cap L^\gamma) (\Omega))
\]
as the ``residual set'' is of small measure, cf. (\ref{u14+}).

\medskip

\noindent Similarly,  we deduce from (\ref{u1++++}), (\ref{u15+}) that
\bFormula{U2}
\{ \vc{m}_\ep \}_{\ep > 0} \ \mbox{is bounded in} \ L^\infty(0,T; L^2 + L^q(\Omega;R^3)), \ q = \frac{2 \gamma}{\gamma + 1}\cdotp
\eF
Moreover, combining (\ref{cutoff}), hypothesis (\ref{u1--}) and (\ref{U1}) we may infer that
\bFormula{U2a}
\left\{ \ep^\alpha \chi_\ep \frac{\vre - \tilde \vre}{\ep^m} \Grad G \right\}_{\ep > 0}
\ \mbox{is bounded in}\ L^\infty(0,T; L^2 + L^1 (\Omega)).
\eF
Finally, by virtue of (\ref{u15+}), (\ref{u16+}),
\bFormula{U3}
\{ \tn{S}(\Grad \vue) \}_{\ep > 0} \ \mbox{is bounded in}\ L^2(0,T; L^2(\Omega; R^{3 \times 3})),
\eF
and
\bFormula{U4}
\{ \vre \vue \otimes \vue \}_{\ep > 0} \ \mbox{is bounded in}\ L^\infty(0,T; L^1(\Omega; R^{3 \times 3})).
\eF

\medskip

\noindent 
Now let us estimate   the pressure perturbation.
Writing
\[
p(\vre) - p(\tilde \vr_{\ep}) - p'(1) (\vre - \tilde \vr_{\ep}) =
p(\vre) - p(\tilde \vr_{\ep}) - p'(\tilde \vr_{\ep}) (\vre - \tilde \vr_{\ep}) +
\Big( p'(\tilde \vr_{\ep}) - p'(1) \Big) (\vre - \tilde \vr_{\ep})
\]
we have
$$
p(\vre) - p(\tilde \vr_{\ep}) - p'(\tilde \vr_{\ep}) (\vre - \tilde \vr_{\ep})
= \Big[ p(\vre) - p(\tilde \vr_{\ep}) - p'(\tilde \vr_{\ep}) (\vre - \tilde \vr_{\ep}) \Big]_{\rm ess} +
\Big[ p(\vre) - p(\tilde \vr_{\ep}) - p'(\tilde \vr_{\ep}) (\vre - \tilde \vr_{\ep}) \Big]_{\rm res}
$$
where, since $p$ is twice continuously differentiable in $(0,\infty)$ and
$\tilde \vre$ satisfies (\ref{odhad}),
\[
\left|
\Big[ p(\vre) - p(\tilde \vr_{\ep}) - p'(\tilde \vr_{\ep}) (\vre - \tilde \vr_{\ep}) \Big]_{\rm ess}(t,x)
\right| \leq c \left[ \vre - \tilde \vr_{\ep} \right]_{\rm ess}^2(t,x) \ \mbox{provided}\ x \in B_{2/\ep^\alpha}.
\]
Similarly, by virtue of (\ref{u1++++}), (\ref{u14+})
\[
{\rm ess} \sup_{t \in (0,T)}  \int_{B_{2/\ep^\alpha} } \left| \Big[ p(\vre) - p(\tilde \vr_{\ep}) - p'(\tilde \vr_{\ep}) (\vre - \tilde \vr_{\ep}) \Big]_{\rm res} \right| \ \dx \leq c \ep^{2m}.
\]
Finally, in accordance with (\ref{odhad}),
\[
| p'(\tilde \vr_{\ep} (x)) - p'(1) | \leq c |\tilde \vre(x) - 1|
\leq c \ep^{2(m-1 - \alpha)} \ \mbox{for} \ x \in B_{2/\ep^\alpha}.
\]
Thus, summing up the previous estimates, we conclude that
\bFormula{estim}
{\rm ess} \sup_{t \in (0,T)} \left\|
\frac{1}{\ep^{2m}} \Big( p(\vre) - p(\tilde \vr_{\ep}) - p'(1)
(\vre - \tilde \vr_{\ep}) \Big)  \right\|_{(L^1 + L^2) (B_{2/\ep^\alpha}; R^3)} \leq c (1 + \ep^{m-2 - 2\alpha} ) .
\eF

\subsubsection{Spatial regularization}

In view of the uniform bounds obtained in the previous section, equations (\ref{aku1r}),
(\ref{aku2r}) can be written in the form
\bFormula{aku1rr}
\ep^m \partial_t S_\ep + \Div \vc{m}_\ep = \ep^\alpha F^1_\ep,
\eF
\bFormula{aku2rr}
\ep^m \partial_t \vc{m}_\ep + p'(1) \Grad S_\ep =
(\ep^{m} +  \ep^{2(m-1-\alpha)}  ) \Div \tn{F}^2_\ep + (\ep^{m -1}  + \ep^\alpha  + \ep^{2(m-1-\alpha)} ) \vc{F}^3_\ep,
\eF
with
\[
\left\{
\begin{array}{c}
\{ F^1_\ep \}_{\ep > 0} \ \mbox{bounded in}\ L^\infty(0,T; L^2 + L^1 (\Omega))\\ \\
\{ \tn{F}^2_\ep \}_{\ep > 0} \ \mbox{bounded in}\ L^2(0,T; L^2 + L^1 (\Omega;R^{3 \times 3})) \\ \\
\{ \vc{F}^3_\ep \}_{\ep > 0} \ \mbox{bounded in}\ L^2(0,T; L^2 + L^1 (\Omega;R^{3 })) .
\end{array} \right\}
\]
Equations (\ref{aku1rr}), (\ref{aku2rr}) are satisfied in the sense of distributions. For future analysis, however, it is more convenient to deal with classical (smooth) solutions. To this end,
we introduce
\[
v_\delta = \kappa_\delta * v,
\]
where $\kappa_\delta = \kappa_\delta(x)$ is a family of regularizing kernels acting in the spatial variable. Accordingly, we may regularize
(\ref{aku1rr}), (\ref{aku2rr}) to obtain
\bFormula{aku1rg}
\ep^m \partial_t S_{\ep,\delta} + \Div \vc{m}_{\ep,\delta} = \ep^\alpha F^1_{\ep,\delta}
\eF
\bFormula{aku2rg}
\ep^m \partial_t \vc{m}_{\ep, \delta} + p'(1) \Grad S_{\ep,\delta} =
(\ep^{m} +  \ep^{2(m-1-\alpha)}  ) \Div \tn{F}^2_{\ep,\delta} + (\ep^{m -1}  + \ep^\alpha  + \ep^{2(m-1-\alpha)} ) \vc{F}^3_{\ep,\delta}
\eF
where
\bFormula{aku3r}
\left\{
\begin{array}{c}
\| {F}^1_{\ep, \delta} \|_{L^2(0,T; W^{k,2}(\Omega)} \leq c(k, \delta), \\ \\ \| \tn{F}^2_{\ep, \delta} \|_{L^2(0,T; W^{k,2}(\Omega; R^{3\times 3}))} \leq c(k, \delta), \\ \\
\| \vc{F}^3_{\ep, \delta} \|_{L^2(0,T; W^{k,2}(\Omega;R^3)} \leq c(k, \delta)
\end{array}
\right\}
\eF
for any $k=0,1,\dots$
uniformly for $\ep \to 0$.

\subsubsection{New formulation of the regularized acoustic wave equation}

Our aim is to rewrite the regularized acoustic equation (\ref{aku1rg}), (\ref{aku2rg}) in terms of the acoustic potential. To this end, we decompose
\[
\vc{m}_{\ep,\delta} = \vc{Y}_{\ep, \delta} + \Grad \Psi_{\ep,\delta},
\]
where $\vc{Y}_{\ep, \delta} = \vc{H}[\vc{m}_{\ep, \delta}]$.

\medskip
\noindent
Introducing  new unknowns
$S_{\ep, \delta}$, $\Psi_{\ep,\delta}$, we may rewrite the acoustic equation
(\ref{aku1rg}), (\ref{aku2rg}) in the form:
\bFormula{a1}
\ep^m \partial_t S_{\ep,\delta} + \Delta \Psi_{\ep,\delta} = \ep^\alpha F^1_{\ep, \delta},
\eF
\bFormula{a2}
\ep^m \partial_t \Psi_{\ep, \delta} + p'(1) S_{\ep,\delta} = (\ep^m + \ep^{2(m-1-\alpha)} ) \Delta^{-1} \Div \Div
[ \tn{F}^2_{\ep, \delta} ] + (\ep^{m-1} + \ep^\alpha +  \ep^{2(m-1-\alpha)} ) \Delta^{-1} \Div \vc{F}^3_{\ep, \delta}.
\eF
Our aim is to use dispersive estimates for (\ref{a1}), (\ref{a2}) to deduce local decay of the acoustic potential. To this end, we present a simple result, the proof of which is a straightforward adaptation of Metcalfe~\cite[Lemma 4.1]{Metca} {(cf. also D'Ancona and Racke~\cite[Example 1.2]{DancRac})} :

\bLemma{w1}
Consider $\varphi \in \DC(R^2)$.
Then
\bFormula{iw1}
\int_{-\infty}^{\infty} \intO{ \left| \varphi (x_h) \exp \left( {\rm i} \sqrt{-\Delta} t \right)[v] \right|^2 } \ \dt \leq c(\varphi) \| v \|^2_{L^2(\Omega)}.
\eF
\eL
\bProof
For a function $w=w(t,x_h, x_3)$, we denote by $\widehat w (\tau, \xi ,k)$ its
Fourier transform in the \emph{space-time} variables, $ \tau \in R$, $\xi \in R^2$, $k \in Z$. Accordingly, by virtue of
Parseval's identity,
$$
\begin{aligned}
&\int_{-\infty}^{\infty} \intO{ \left| \varphi (x_h) \exp \left( {\rm i} \sqrt{-\Delta} t \right)[v] \right|^2 } \ \dt  \\
& \quad = c \sum_{k \in Z} \int_{-\infty}^{\infty} \int_{R^2} \left| \int_{R^2} \widehat \varphi(\xi - \eta)
\delta(\tau - \sqrt{|\eta|^2 + k^2}) \widehat v(\eta, k) \ {\rm d}\eta \right|^2 {\rm d} \xi \ {\rm d}\tau \\
&\quad= c \sum_{k \in Z} \int_{-\infty}^{\infty} \int_{R^2} \left| \int_{\{ \tau = \sqrt{ |\eta|^2 + k^2} \} } \widehat \varphi(\xi - \eta)
 \widehat v(\eta, k) \ {\rm d}S_\eta \right|^2 {\rm d} \xi \ {\rm d}\tau.\end{aligned}
$$
Furthermore, by the Cauchy-Schwartz inequality,
\[
\begin{aligned}
&\sum_{k \in Z} \int_{-\infty}^{\infty} \int_{R^2} \left| \int_{\{ \tau = \sqrt{ |\eta|^2 + k^2} \} } \widehat \varphi(\xi - \eta)
 \widehat v(\eta, k) \ {\rm d}S_\eta \right|^2 {\rm d} \xi \ {\rm d}\tau
\\
&\quad\leq
\sum_{k \in Z} \int_{-\infty}^{\infty} \int_{R^2} \left( \int_{\{ \tau = \sqrt{ |\eta|^2 + k^2} \} } |\widehat \varphi(\xi - \eta)| \ {\rm d}S_\eta \right) \left(  \int_{\{ \tau = \sqrt{ |\eta|^2 + k^2} \} } |\widehat \varphi(\xi - \eta)|
| \widehat v(\eta, k)|^2 \ {\rm d}S_\eta \right) {\rm d} \xi \ {\rm d}\tau
\\
&\quad\leq c(\varphi) \sum_{k \in Z} \int_{R^2} \int_{-\infty}^\infty \int_{\{ \tau = \sqrt{ |\eta|^2 + k^2} \} } |\widehat \varphi(\xi - \eta)|
| \widehat v(\eta, k)|^2 \ {\rm d}S_\eta \ {\rm d}\tau \ {\rm d}\xi
\\
&\quad\leq c(\varphi) \sum_{k \in Z} \int_{R^2}  \int_{R^2} |\widehat \varphi(\xi - \eta)|
| \widehat v(\eta, k)|^2 \ {\rm d}\eta \ {\rm d}\xi \leq c(\varphi)
\| v \|_{L^2(\Omega)}^2.
\end{aligned}
\]
That proves Lemma~\ref{Lw1}.

\qed

\bigskip
\noindent
The gradient of the acoustic potential $\Psi_{\ep, \delta}$ can be written by means of Duhamel's formula:
$$
\begin{aligned}
\Grad \Psi_{\ep,\delta} &= \frac{1}{2} \exp \left({\rm i} \sqrt{-\Delta} \frac{t}{\ep^m}
\right) \left[ \Grad \Psi_{0,\ep, \delta} + {\rm i} \frac{1}{\sqrt{-\Delta}} \Grad [S_{0,\ep, \delta}]
\right]
\\
& {} \quad + \frac{1}{2} \exp \left(-{\rm i} \sqrt{-\Delta} \frac{t}{\ep^m} \right) \left[
\Grad \Psi_{0,\ep, \delta} - {\rm i} \frac{1}{\sqrt{-\Delta}} \Grad [S_{0,\ep, \delta}] \right]
\\
& {} \quad 
+ \frac{\ep^{\alpha -m }}{2} \int_0^t \left( \exp \left({\rm i} \sqrt{- \Delta} \frac{t - s}{\ep^m}
\right) - \exp \left(-{\rm i} \sqrt{-\Delta} \frac{t - s}{\ep^m}
\right)\right)\left[ \frac{{\rm i}}{\sqrt{-\Delta}} \Grad  {F}^1_{\ep, \delta} \right] \ {\rm d}s \\
& {} \quad +
\frac{1 + \ep^{m - 2 - 2\alpha}}{2} \int_0^t \left( \exp \left({\rm i} \sqrt{- \Delta} \frac{t - s}{\ep^m}
\right) + \exp \left(-{\rm i} \sqrt{-\Delta} \frac{t - s}{\ep^m}
\right)\right)\left[ \Grad \Delta^{-1} \Div \Div \tn{F}^2_{\ep, \delta} \right] \ {\rm d}s
\\
& {} \quad +
\frac{\ep^{-1} + \ep^{\alpha - m} + \ep^{m - 2 - 2\alpha}}{2} \!\int_0^t \left(\! \exp \left({\rm i} \sqrt{- \Delta} \frac{t - s}{\ep^m}
\right) \!\! +\!\exp \left(-{\rm i} \sqrt{-\Delta} \frac{t - s}{\ep^m}
\right)\!\!\right)\left[ \Grad \Delta^{-1} \Div \vc{F}^3_{\ep, \delta} \right] {\rm d}s
\end{aligned}
$$
where, for the sake of simplicity, we have set $p'(1) = 1$.
\medskip
\noindent
Now, in accordance with Lemma \ref{Lw1},
\bFormula{od1}
\int_0^T \int_K \left| \exp\left( {\rm i} \sqrt{-\Delta} \frac{t}{\ep^m} \right) [v]
\right|^2 {\rm d}x \ {\rm d}t \leq \ep^m \int_0^\infty \int_K \left| \exp\left( {\rm i} \sqrt{-\Delta} {t} \right) [v]
\right|^2 {\rm d}x \ {\rm d}t \leq \ep^m c \| v \|^2_{L^2(\Omega)},
\eF
and, similarly,
\bFormula{od2}
\begin{aligned}
&\int_0^T \int_K \left| \int_0^t \exp \left( {\rm i} \sqrt{-\Delta} \frac{t-s}{\ep^m} \right)  [g(s)] \ {\rm d}s \right|^2 \ {\rm d}x \ {\rm d}t \\
&\quad \leq T \int_0^T \int_0^T \int_K \left| \exp \left( {\rm i} \sqrt{-\Delta} \frac{t-s}{\ep^m} \right)  [g(s)] \ \right|^2 \ {\rm d}x \ {\rm d}t \ {\rm d}s
\\
&\quad\leq c T  \ep^m   \int_0^T \left\| \exp\left( - {\rm i} \frac{s}{\ep} \right) [g(s)] \right\|^2_{L^2(\Omega)} = \ep^m \| g \|_{L^2((0,T) \times \Omega)}^2
\end{aligned}
\eF
for any compact $K \subset \Omega$.
\medskip
\noindent
Combining (\ref{od1}), (\ref{od2}) with the uniform bounds (\ref{aku3r}) and hypotheses~(\ref{u1+}),~(\ref{u1++}) we infer that
\bFormula{od3}
\int_0^T
\| \Grad \Psi_{\ep,\delta} \|_{L^{2}(K;R^3)}^2 \ \dt \leq \ep^{2\beta} c(\delta, K,T)
\ \mbox{for any compact}\ K \subset \Omega
\eF
uniformly for $\ep \to 0$, where $\beta > 1$ provided $\alpha$ satisfies
(\ref{akuhyp}). Thus the effect of acoustic waves becomes negligible in the
limit $\ep \to 0$. { Accordingly, the information about the asymptotic behavior is
provided by the solenoidal component of the velocity field analyzed in the following section.}

\subsection{Solenoidal part}

\label{c}

In order to control the solenoidal component of the velocity field, we write the momentum
equation (\ref{aku2r}) in the form:
$$
\begin{aligned}\ep \partial_t \vc{m}_{\ep,\delta}  +   \vc{b} \times \vc{m}_{\ep,\delta}
 = \ep \Div \Big[ \chi_\ep \tn{S} (\Grad \vue) - \chi_\ep (\vre \vue \times \vue) \Big]_\delta
- \ep \Big[ \Big( \tn{S} (\Grad \vue) - \vre \vue \times \vue \Big) \cdot \Grad \chi_\ep
\Big]_\delta
\\
+ \ep^{m-1} \Big[ \chi_\ep \frac{ \vre - \tilde \vr_{\ep} }{\ep^m}  \Grad G \Big]_\delta
+ \ep^{1 - 2m} \Big( \Grad [ \chi_\ep (p(\tilde \vre) - p(\vre))]_\delta -
[ \Grad \chi_\ep  ( p(\tilde \vre) - p(\vre) ) ]_\delta \Big)
\end{aligned}
$$
in other words
\bFormula{aku2rm}
\ep \partial_t \vc{m}_{\ep,\delta}  +   \vc{b} \times \vc{m}_{\ep,\delta}= (\ep + \ep^{m-1-\alpha}) \vc{Q}_{\ep, \delta} + \ep^{1 - 2m} \Big( \Grad [ \chi_\ep (p(\tilde \vre) - p(\vre))]_\delta -
[ \Grad \chi_\ep  ( p(\tilde \vre) - p(\vre) ) ]_\delta \Big),
\eF
where
\[
\{ \vc{Q}_{\ep, \delta} \}_{\ep > 0} \ \mbox{is bounded in}\ L^2(0,T; W^{k,2}(\Omega;R^3))
\ \mbox{for any fixed}\ k , \ \delta > 0.
\]

\subsubsection{Compactness of vertical averages}
Taking the vertical average of equation (\ref{aku2rm}) we obtain
\bFormula{vaeq}
\begin{aligned}\ep \partial_t \va{ \vc{m}_{\ep,\delta} }  +   \vc{b} \times \va{ \vc{m}_{\ep,\delta}}
&=  (\ep + \ep^{m-1-\alpha}) \va{ \vc{Q}_{\ep, \delta} }\\
&\quad + \ep^{1 - 2m} \Big( \Grad \va{[ \chi_\ep (p(\tilde \vre) - p(\vre))]_\delta} -
\va{[ \Grad \chi_\ep  ( p(\tilde \vre) - p(\vre) ) ]_\delta} \Big).
\end{aligned}
\eF
Recalling
\[
\vc{m}_{\ep, \delta} = \vc{Y}_{\ep, \delta} + \Grad \Psi_{\ep, \delta},\
\vc{Y}_{\ep, \delta} \equiv \vc{H}[ \vc{m}_{\ep, \delta} ]
\]
we check easily that
\[
\vc{b} \times \va{ \vc{Y}_{\ep, \delta} } =
\left[ \begin{array}{c} - \va{Y^2_{\ep,\delta}} \\ \va{Y^1_{\ep,\delta}} \\ 0 \end{array} \right]
\]
is an exact gradient as $\vc{Y}_{\ep, \delta}$ is solenoidal.

\medskip
\noindent
Consequently, testing equation (\ref{vaeq}) on $\varphi \in \DC(\Omega; R^3)$, $\Div \varphi = 0$ we get
\[
\partial_t \intO{ \va{\vc{m}_{\ep, \delta}} \cdot \varphi } =
\intO{ \va{ \vc{Q}_{\ep, \delta} } \cdot \varphi } - \frac{1}{\ep} \intO{
\va{ \vc{b} \times \Grad \Psi_{\ep,\delta} } \cdot \varphi }
\]
provided $\ep$ is small enough so that $\chi_\ep |_{{\rm supp} \varphi} \equiv 1$.
Thus we may use (\ref{od3}) to conclude that
\bFormula{i5}
\va{\vc{Y}_{\ep, \delta}} \to \vc{U}_\delta \ \mbox{strongly in}\ L^2((0,T) \times K; R^2)
\ \mbox{for any compact} \ K \subset \Omega \ \mbox{and any fixed}\ \delta > 0.
\eF
We note that this step depends essentially on (\ref{akuhyp}) that requires the rather
strong assumption $m > 10$.

\subsubsection{Oscillations}
We write any function $v$ in the form
\begin{equation}\label{[v]}
v(x) = \va{v}(x_h) + \Bl v \Br (x).
\end{equation}
Since $ \Bl v \Br = v - \va{v}$ has zero vertical mean, it can be written in the form
\[
\Bl v \Br (x) \equiv \partial_{x_3} I[v] \quad\mbox{with}\quad\int_{{\cal
T}^1} I[v](x) \ {\rm d}x_3 = 0.
\]
Moreover, we define
\[
\omega^{i,j}_{\ep,\delta} = \partial_{x_i} {Y}_{\ep,\delta}^j - \partial_{x_j}
{Y}^i_{\ep, \delta} = \partial_{x_i} m^j_{\ep, \delta} - \partial_{x_j} m^i_{\ep, \delta}.
\]
Going back to equation (\ref{aku2rm}), we deduce that
\bFormula{i6}
\ep \partial_t \omega^{1,2}_{\ep, \delta} + {\rm div}_h [\vc{Y}_{\ep,\delta}]_h
= \ep \Big( \partial_{x_1} Q^2_{\ep, \delta} - \partial_{x_2} Q^1_{\ep, \delta} \Big)
- \Delta_h \Psi_{\ep, \delta},
\eF
\bFormula{i6A}
\ep \partial_t \omega^{1,3}_{\ep, \delta} + \partial_{x_3} Y^2_{\ep, \delta} =
\ep \Big( \partial_{x_1} Q^3_{\ep, \delta} - \partial_{x_3} Q^1_{\ep, \delta} \Big)
- \partial^2_{x_3,x_2} \Psi_{\ep, \delta},
\eF
and, finally,
\bFormula{i6B}
\ep \partial_t \omega^{2,3}_{\ep, \delta} - \partial_{x_3} Y^1_{\ep, \delta} =
\ep \Big( \partial_{x_2} Q^3_{\ep, \delta} - \partial_{x_3} Q^2_{\ep, \delta} \Big)
+ \partial^2_{x_3,x_1} \Psi_{\ep, \delta}
\eF
for all $t \in (0,T)$, $x \in B_{1/\ep^\alpha}$.

\subsubsection{Analysis of the convective term}

Following step by step the analysis performed in \cite[Section 3]{GallSR} we observe that
the only problematic component of the convective term reads
\[
\vc{H} \left[
\int_{ {\cal T}^1} \left[ \Div ( \vc{Y}_{\ep, \delta} \otimes \vc{Y}_{\ep, \delta}  ) \right]_h \ {\rm d}x_3 \right].
\]

$\bullet$ {\bf Step 1:}

\medskip

Since $\vc{Y}_{\ep, \delta}$ is solenoidal,
we have
\[
\Div (\vc{Y}_{\ep, \delta} \otimes \vc{Y}_{\ep, \delta} ) = \frac{1}{2} \Grad |\vc{Y}_{\ep, \delta} |^2 - \vc{Y}_{\ep, \delta} \times ({\bf curl}[\vc{Y}_{\ep, \delta} ]).
\]
As the former term is a gradient, we concentrate on the latter.

\bigskip

$\bullet$ {\bf Step 2:}

\medskip

Write
\[
\begin{aligned}
\vc{Y}_{\ep, \delta} \times  {\bf curl}[ \vc{Y}_{\ep, \delta} ] &=
\va{\vc{Y}_{\ep, \delta}} \times {\bf curl} \va{\vc{Y}_{\ep, \delta}} + \partial_{x_3} \Big(
\va{\vc{Y}_{\ep, \delta}} \times  {\bf curl} \ I [{\vc{Y}}_{\ep, \delta}] + I[\vc{Y}_{\ep, \delta}] \times  {\bf curl} \va{\vc{Y}_{\ep, \delta}}
\Big)
\\
& {} \qquad +\partial_{x_3} I[\vc{Y}_{\ep, \delta}] \times \partial_{x_3} {\bf curl}   (I[\vc{Y}_{\ep,\delta}]),
\end{aligned}\]
where the term in the brackets has zero vertical mean, while, in accordance with (\ref{i5}), the first term is pre-compact. Finally, we have
\[\begin{aligned}
&
\Big[
\partial_{x_3} I[\vc{Y}_{\ep, \delta}] \times \partial_{x_3} {\bf curl}   (I[\vc{Y}_{\ep,\delta}]) \Big]^j
\\
& = \partial_{x_3} I[Y^i_{\ep, \delta}]
\partial_{x_3} \Big( \partial_{x_i} I[Y^j_{\ep, \delta}] - \partial_{x_j} I[Y^i_{\ep, \delta}] \Big) =
\partial_{x_3} I[Y^i_{\ep, \delta}]
\partial_{x_3} I[ \omega^{i,j}_{\ep, \delta}]
, \ j=1,2,3.
\end{aligned}\]

\bigskip

$\bullet$ {\bf Step 3:}

\medskip

In accordance with (\ref{i6} - \ref{i6B}), we get
\[
\ep \partial_t (\partial_{x_3} I[\omega^{1,3}_{\ep, \delta}] ) + \partial^2_{x_3}
I[Y^2_{\ep, \delta}] = \ep \Big( \partial_{x_1}( Q^3_{\ep, \delta} - \va{ Q^3_{\ep, \delta} }) - \partial_{x_3} Q^1_{\ep, \delta} \Big) - \partial^2_{x_3,x_2} \Psi_{\ep, \delta},
\]
and
\[
\ep \partial_t (\partial_{x_3} I[\omega^{2,3}_{\ep, \delta}] ) - \partial^2_{x_3}
I[Y^1_{\ep, \delta}] = \ep \Big( \partial_{x_2}( Q^3_{\ep, \delta} - \va{ Q^3_{\ep, \delta} }) - \partial_{x_3} Q^2_{\ep, \delta} \Big) - \partial^2_{x_3,x_1} \Psi_{\ep, \delta}
\]
at least for $x \in B_{1/\ep^\alpha}$.
Next, compute
\[\begin{aligned}
\left[
\partial_{x_3} I[\vc{Y}_{\ep, \delta}] \times \partial_{x_3} {\bf curl}   (I[\vc{Y}_{\ep,\delta}]) \right]^1 = \partial_{x_3} I[Y^2_{\ep, \delta}] \partial_{x_3}
I[\omega^{2,1}_{\ep, \delta}] + \partial_{x_3} I[Y^3_{\ep, \delta}] \partial_{x_3}
I[\omega^{3,1}_{\ep, \delta}]
\\
=
\partial_{x_3} \Big( \partial_{x_3} I[Y^2_{\ep, \delta}]
I[\omega^{2,1}_{\ep, \delta}] \Big) - \partial^2_{x_3} I[Y^2_{\ep, \delta}]I[\omega^{2,1}_{\ep, \delta}] - I [{\rm div}_h [\vc{Y}_{\ep, \delta}]_h ]
\partial_{x_3} I[\omega^{3,1}_{\ep, \delta} ],
\end{aligned}\]
where, furthermore,
\[
\begin{aligned}
\partial^2_{x_3} I[Y^2_{\ep, \delta}]I[\omega^{2,1}_{\ep, \delta}] + I [{\rm div}_h [\vc{Y}_{\ep, \delta}]_h ]
\partial_{x_3} I[\omega^{3,1}_{\ep, \delta} ] = - \ep \partial_t \Big( \partial_{x_3}
I[\omega^{1,3}_{\ep, \delta} ] \Big) I[\omega^{2,1}_{\ep, \delta}]
\\
+ \ep \Big( \partial_{x_1}( Q^3_{\ep, \delta} - \va{ Q^3 }_{\ep, \delta}) - \partial_{x_3} Q^1_{\ep, \delta} \Big) I[\omega^{2,1}_{\ep, \delta}]- \partial^2_{x_3,x_2} \Psi_{\ep, \delta} I[\omega^{2,1}_{\ep, \delta}] - \ep \partial_t
\Big( I[ \omega^{1,2}_{\ep, \delta}] \Big) \partial_{x_3} I[ \omega^{3,1}_{\ep, \delta}]
\\
+ \ep I \Big[ \partial_{x_1} Q^2_{\ep, \delta} - \partial_{x_2} Q^1_{\ep, \delta} \Big]
\partial_{x_3} I[\omega^{3,1}_{\ep,\delta} ] - I [ \Delta_h \Psi_{\ep, \delta} ] I[\omega^{3,1}_{\ep,\delta} ]
\\
= \ep \partial_t \Big( \partial_{x_3} I[\omega^{1,3}_{\ep, \delta} ] I[\omega^{1,2}_{\ep, \delta}] \Big) + \ep \Big( \partial_{x_1}( Q^3_{\ep, \delta} - \va{ Q^3 }_{\ep, \delta}) - \partial_{x_3} Q^1_{\ep, \delta} \Big) I[\omega^{2,1}_{\ep, \delta}]- \partial^2_{x_3,x_2} \Psi_{\ep, \delta} I[\omega^{2,1}_{\ep, \delta}]
\\
+ \ep I\Big[ \partial_{x_1} Q^2_{\ep, \delta} - \partial_{x_2} Q^1_{\ep, \delta} \Big]
\partial_{x_3} I[\omega^{3,1}_{\ep,\delta} ] - I [ \Delta_h \Psi_{\ep, \delta} ] I[\omega^{3,1}_{\ep,\delta} ].
\end{aligned}
\]
for $t \in (0,T)$, $x \in B_{1/\ep^\alpha}$.
Treating the second component in a similar fashion, we conclude that
\bFormula{converg}
\frac{1}{|{\cal T}^1|} \int_0^T \intO{
 \left[ \Div ( \vc{Y}_{\ep, \delta} \otimes \vc{Y}_{\ep, \delta}  ) \right] \cdot \varphi } \ \dt \to  \int_0^T \int_{R^2} \Div ( \vc{U}_{\delta} \otimes \vc{U}_{\delta}  ) \cdot \phi
 \ {\rm d}x_h \ \dt \ \mbox{as}\ \ep \to 0
\eF
for any fixed $\delta > 0$, and for any $\varphi = [\phi(t,x_h), 0]$,
$\phi \in \DC([0,T) \times R^2;R^2)$, ${\rm div}_h \phi = 0$.

\subsection{Convergence - conclusion}

In order to complete the proof of Theorem \ref{Tm1}, write
\[
\begin{aligned}
\vre \vue \otimes \vue
= (\vre - 1) \vue \otimes \vue +  ( \vu_{\ep, \delta} - \vue ) \otimes \vu_{\ep, \delta} +  \vu_{\ep, \delta}  \otimes (\vu_{\ep, \delta} - \vue)
\\
{} + [(1 - \vre) \vue]_{\delta}  \otimes \vu_{\ep, \delta} +
\vu_{\ep,\delta} \otimes [(1 - \vre) \vue]_{\delta}
 \vc{m}_{\ep, \delta} \otimes \vc{m}_{\ep, \delta}
\end{aligned}
\]
for $t \in (0,T)$ and $x \in K \subset \Omega$, $K$ compact.

\medskip
\noindent
Now, it is easy to observe that
\[
(\vre - 1) \vue \otimes \vue, \ [(1 - \vre) \vue]_{\delta}  \otimes \vu_{\ep, \delta} \to 0 \ \mbox{in}\ L^1((0,T) \times K) \ \mbox{for}\ \ep \to 0
\]
while
\[
\| \vu_{\ep, \delta}   - \vue   \|_{L^2(K;R^3)} \leq c \delta \| \vue \|_{W^{1,2}(K;R^3)}
\ \mbox{uniformly for} \ \ep > 0.
\]
Consequently, testing the momentum equation (\ref{i2}) on a compactly supported solenoidal function, the convective term $\vre \vue \otimes \vue$ can be replaced by
$\vc{m}_{\ep, \delta} \otimes \vc{m}_{\ep, \delta}$ handled in detail in Sections
\ref{a}, \ref{c}. This completes the proof of Theorem \ref{Tm1}. 

\qed

\medskip
\noindent
Note that, since the limit problem (\ref{NS1}), (\ref{NS2}) admits a \emph{unique} solution for any square integrable initial data, there is no need to consider subsequences.

\section{Isotropic scaling: Proof of Theorem~\ref{Tm2}}

\label{proof2}

We consider the situation where the Rossby and Mach numbers have the same scaling.
Accordingly, we have $m=1$, together with the uniform bounds derived in Paragraph~\ref{u}. We recall that the static solution~$\tilde \vre = \tilde \vr (|x_h|) $ is now independent of $\ep$ and we set
\[
r_\ep = \frac{\vre - \tilde \vr}{\ep}\cdotp
\]

\subsection{Preliminary results}
%
%\label{u1+++1} {\rm ess} \sup_{t \in (0,T)} \left\|  \left[ \frac{\vre - \tilde \vr_{\ep} }{\ep^m} \right]_{\rm ess} \right\|_{L^2(B_{r})} \leq c(r) & \\ 
%\label{u1++++1}{\rm ess} \sup_{t \in (0,T)} \int_{B_{r}} [\vre]^\gamma_{\rm res} \ \dx \leq \ep^{2m} c(r)& \\ 
%\label {u14+1}
%{\rm ess} \sup_{t \in (0,T)} \int_{B_{r}}  1_{\rm res} \ \dx \leq \ep^{2m} c
%(r).&
%\end{eqnarray}
%Finally, combining (\ref{u14++} - \ref{u15++})

In accordance with the uniform bounds established in~(\ref{u1+++1}-\ref{u14+1}) and~(\ref{u14++}-\ref{u15++}),
we may assume (up to taking subsequences) that
\begin{equation} \label{weakconvrho}
 \vre \rightarrow  \tilde \vr \:  \mbox{  in } L^\infty(0,T; L^\gamma(K)), \qquad r_\eps \rightarrow r \: \mbox{ weakly-$(^*)$ in } L^\infty(0,T; L^2(K))
 \end{equation}
for any compact $K \subset \Omega$, recalling that $\gamma > 3$. Moreover,
\begin{equation} \label{weakconvu}
\vue \rightarrow \vc{U} \:  \mbox{ weakly in } \: L^2(0,T; W^{1,2}(\Omega; R^3)).
\end{equation}
The main point is to derive the equations satisfied by $r$ and $\vc{U}$. Similarly to
Section \ref{a}, we rewrite the Navier-Stokes system in the following form:
\begin{equation} \label{NSbis}
\left\{
\begin{aligned}
& \eps \pa_t r_\eps + \Div (\vre \vue)   = 0, \\
& \eps \pa_t (\vre \vue) + \eps \Div \left( \vre \vue \otimes \vue \right) + \Grad \left( p'(\tilde \vr)  r_\eps \right) + \vc{b} \times (\vre \vue)  \\
& = - \eps \Div \tn{S}(\Grad \vue) - \frac{1}{\eps} \Grad \left( p(\vre) - p(\tilde \vr) -  p'(\tilde \vr) (\vre - \tilde \vr)\right) +  r_\eps \Grad G .
\end{aligned}
\right.
\end{equation}
The bounds~(\ref{u14++}-\ref{u14+1})   allow to pass to the limit in the previous system, in the sense of distributions, to derive the following diagnostic equations:
$$\Div(\tilde \vr \vc{U}) = 0, \quad  \Grad(p'(\tilde \vr) \, r) + \vc{b} \times \tilde \vr \vc{U} =  r \Grad G.   $$
Moreover, introducing
\[
R \equiv  P'(\tilde \vr)  r,
\]
one can write the previous equations as
\begin{equation*} \Div(\tilde \vr \vc{U}) = 0, \quad   \Grad R  + \vc{b} \times   \vc{U} = 0.
\end{equation*}
From there, one has easily that
\begin{equation}\label{A}
R = R(x_h), \quad \vc{U} = [\vc{U}_h(x_h),0]
\end{equation}
together with
\begin{equation}\label{B}
\Divh \left( \tilde \vr \vc{U}_h\right)  =0, \quad  \na_h R  +
\vc{U}_h^\perp = 0.
\end{equation}
Note that in addition to the
constraint $\Divh \left( \tilde \vr \vc{U}_h\right)  = 0$, the last
equation implies that~$\Divh \vc{U}_h = 0$, and~$\Grad \tilde \vr \cdot \vc{U}_h = 0$, therefore
$R = R(|x_h|)$ is radially symmetric.

\medskip
\noindent
These
constraints do not determine the time evolution of $r$ and $\vc{U}_h$.
Our aim is to show that, in accordance with the conclusion of Theorem \ref{Tm2},
the function $r$ satisfies equation (\ref{limiteq}).

\medskip
\noindent
{The rest of the paper is  devoted to the proof of (\ref{limiteq}).
To this end, we perform the
analysis of the high  frequency waves generated by the linear part of system
(\ref{NSbis}).} Indeed, one can  recast   \eqref{NSbis} in the general form
$$
 \eps \pa_t {\cal U}_\eps  + L[ {\cal U}_\eps] = \eps {\cal N}_\eps + \eps{\cal F}_\eps
 $$
where  ${\cal U}_\eps \equiv [r_\eps,\vc V_\eps]$, with $\vc V_\eps \equiv \vre
\vue$,  while
$$L [{\cal U}_\eps] \equiv  \left[ \Div \vc V_\eps ,  \,  \tilde \vr \, \Grad(P'(\tilde \vr) \, r_\eps) + \vc{b}
\times \vc V_\eps \right]$$
stands for the linear part of the equation, and
 \begin{equation}\label{defNep}
{\cal N}_\eps = [0, \vc N_\ep ],\quad \vc N_\ep=  -\Div \left( \vre \vue
\otimes \vue \right) .
 \end{equation}
Furthermore,  ${\cal F}_\eps$ includes the remaining terms:
 \begin{equation}\label{defSep}
{\cal F}_\eps =  (0, \vc{F}_\ep),\quad \vc{F}_\ep= - \Div \tn{S}(\Grad \vue)
- \frac{1}{\eps^2} \Grad \Big( p(\vre) - p(\tilde \vr) -
p'(\tilde \vr) (\vre - \tilde \vr)\Big) .
  \end{equation}
It can be checked that $L$ defines a skew-symmetric operator with respect to  the scalar product
 $$ \left( {\cal U} ,  {\cal U}' \right)  \: \equiv \:  \int_{\Omega} \Big(  r r' P'(\tilde \vr) \: + \: \vc V \cdot\vc  V' \, \tilde \vr^{-1} \Big) \: 
 {\rm d}x;  $$
whence it generates time oscillations with frequency
$\ep^{-1}$, whose nonlinear interaction is not obvious. In fact, {the
the principal difficulty of this part of the paper is the fact that the
spectral properties of the operator $L$ are unclear, and thus, in
contrast with the previous part, employing   dispersive effects
seems difficult if not impossible. We shall therefore rely on
local methods and use as much as possible the structure of the
equations, to control the asymptotic behavior of the convective
term.}

\subsection{Handling the convective term} \label{subsecconvective}

The goal of this section is to better understand the behavior of
$\vc N_\ep$, and more precisely of
\begin{equation}\label{N}
\vc{\widetilde N}_{\ep} \equiv \Divh \left( \frac{1}{\tilde \vr} \va{\vc
V_{\ep,h} \otimes \vc V_{\ep,h}} \right),\quad\mbox{where}\quad\vc
V_{\ep,h}=\vre \vc u_{\ep,h}.
\end{equation}
Later on, we shall compare this term with the average with respect
to $x_3$ of
$$
- \vc N_{\ep,h}=\Div\left(\vre\vue\otimes\vu_{\ep,h}\right),
$$
which reads
\begin{equation}\label{Nh}
\va {- \vc N_{\ep,h}}=\Divh \left(\va{\vre
\vu_{\ep,h}\otimes\vu_{\ep,h}}\right).
\end{equation}
These are two main steps in deriving the limit equation
\eqref{limiteq}.
We start by showing  that
\[
\left(
\vc {\tilde N}_\ep | \na^\perp_h \psi \right) \to 0, \quad \ep \to 0,
\]
for all $\psi \in \DC ( (0,T) \times R^2)$ such that $\na_h \tilde \vr \cdot  \na_h^{\perp} \psi = 0$. Here and hereafter, the symbol $( \cdot | \cdot )$ denotes
the standard duality pairing.
We follow a general  approach initiated by Lions and Masmoudi \cite{LIMA1}, \cite{LIMA6}
in the framework of the incompressible limit, and later adapted in \cite{GallSR} to the case of rotating fluids. This approach, used also in the last step of the proof of Theorem~\ref{Tm1}, is reminiscent of compensated compactness arguments. It relies on various cancellations, obtained by multiple use of the structure of \eqref{NSbis} which roughly speaking reads
 \begin{equation*}
  \left\{
\begin{aligned}
& \eps \pa_t r_\eps + \Div (\vre \vue)   = 0, \\
& \eps \pa_t (\vre \vue) +  \tilde \vr \, \Grad \left( P'(\tilde \vr)  r_\eps \right) + \vc{b} \times (\vre \vue)  =  O(\eps)
\end{aligned}
\right.
\end{equation*}
and will provide some compactness for appropriate  combinations of $r_\ep, \vre, \vue$.

\medskip
\noindent
Similarly to Section \ref{a},
the complete treatment of $\vc{\widetilde N}_\ep$ involves spatial regularizations:
\[
\vu_{\ep,\delta} = \kappa_\delta *\vu_\ep ,
\quad \vc V_{\ep,\delta} \: = \: \kappa_\delta * (\vre \vue), \quad
r_{\ep,\delta} \: = \: \kappa_\delta * r_\ep,
\]
where $\kappa_\delta = \kappa_\delta(x)$ is a family of regularizing kernels.
Note that the quantities
$\vc V_{\ep,\delta} $ and~$ r_{\ep,\delta}$ are bounded
{in~$L^\infty(0,T;L^1 (K) )$ } for any compact $K \subset \Omega$,
{and
$\vu_{\ep,\delta}$ {in~$L^2(0,T; W^{1,2}(\Omega;R^3))$}
},
uniformly in~$\delta,\ep$.
We claim the following result.

\bProposition{rop2} The fields $\vu_{\ep,\delta}$, $\vc
V_{\ep,\delta}$ and $r_{\ep,\delta}$ satisfy the following
properties:
\begin{equation}\label{estimateV1}
  {\vc V_{\eps , \delta}} =   \eps \vc t_{\ep,\delta} ^1 + \vc t_{\ep,\delta}
  ^2,
\quad \mbox{and} \quad
  \curl \left( \frac1{\tilde \vr}\  \vc V_{\eps , \delta} \right ) = \eps \vc T_{\ep,\delta} ^1 +\vc T_{\ep,\delta} ^2
\end{equation}
where
{\begin{equation}\label{estimatecurlV2}
\begin{array}{r}
\| \vc t_{\ep,\delta} ^1 \|_{L^2(0,T ; W^{k,2}(K;R^3))} + \| \vc
T_{\ep,\delta} ^1 \|_{L^2(0,T ; W^{k,2} (K;R^3))} \:\le \: c (\delta),
 \\ \\
\|\vc t_{\ep,\delta} ^2 \|_{L^2(0,T ; W^{1,2}( K;R^3))} + \| \vc
T_{\ep,\delta} ^2 \|_{L^2((0,T) \times K;R^3)} \le c \:
\end{array}
 \end{equation}}
for any compact $K \subset \Omega$ and $k=0,1,\dots$
Moreover,
\begin{equation}\label{delta-cv}
\begin{aligned} \sup_{\ep > 0}
\| r_{\ep} - r_{\ep,\delta}\|_{L^\infty( 0,T ; W^{s,2} (K))} \to 0 \hbox{ as } \delta
\to 0\ \hbox{ for all } s<0,\\
\sup_{\ep > 0}
\|\vc u_\eps-\vc u_{\ep,\delta}\|_{L^2(0,T; W^{s,2}(K;R^3))} \to 0
\hbox{ as } \delta \to 0
\  \hbox{ for all } s<1,\\
\sup_{\ep > 0}
\|\vc V_{\ep} - \vc V_{\ep,\delta}\|_{L^2(0,T; W^{s,2} (K;R^3))} \to
0 \hbox{ as } \delta \to 0 \  \hbox{
for all } s<-1/2
\end{aligned}
\end{equation}
for any compact $K \subset \Omega$.
Finally
the following approximate wave equation is satisfied:
\begin{align} \label{NSbisreg1}
& \eps \pa_t r_{\eps,\delta} + \Div \vc V_{\ep,\delta}   = 0, \\
\label{NSbisreg2} & \eps \pa_t \vc V_{\ep,\delta} +  \tilde
\vr\Grad \left( P'(\tilde \vr) \,  r_{\eps,\delta}  \right) +
\vc{b} \times \vc V_{\ep,\delta}  = \eps  \vc F^1_{\eps, \delta}  +
\vc F^2_{\eps, \delta},
\end{align}
where the source terms are smooth and satisfy
\begin{equation}\label{bddS1epsdelta}
\sup_{\ep > 0} {\| \vc F^1_{\eps, \delta} \|_{L^2(0,T ; W^{k,2}(K;R^3))} }\: \le \:
c(\delta),
\end{equation}
\begin{equation}\label{bddS2epsdelta}
\sup_{\ep > 0} \left(
 {\| \vc F^2_{\eps, \delta} \|_{L^2((0,T) \times K;R^3)} \: + \:
 \left\| \curl \bigl ( \frac 1 {\tilde \vr} \vc F^2_{\eps, \delta}  \bigr)\right \|_{L^2((0,T) \times K;R^3)} } \right) \rightarrow 0
¬†\: \mbox{ as $\delta \rightarrow 0$  }
\end{equation}
for any compact $K \subset \Omega$, $k=0,1,\dots$
\eP
We postpone the proof of this proposition to Section \ref{subsecregular}.
This regularization process allows to  establish    the following convergence result.
 \bProposition{rop1}
The nonlinear quantity
\begin{equation}\label{Nhd}
\vc{\widetilde N}_{\ep,\delta} \: = \: \Divh \left( \frac{1}{\tilde
\vr} \va{\vc V_{\ep,\delta,h} \otimes \vc V_{\ep,\delta,h}}
\right)
\end{equation}
satisfies
$$\displaystyle  \lim_{\delta \rightarrow 0} \left( \limsup_{\eps \rightarrow 0}  \: \left|
\left(  \vc{\widetilde N}_{\ep,\delta} | \na^\perp_h \psi \right) \right| \right) \: = \: 0, \:
$$
 for all $\psi \in \DC ( (0,T) \times R^2)$ such that $\na_h \tilde \vr \cdot  \na_h^{\perp} \psi = 0$.
\eP
The rest of this part is devoted to the proof of Proposition \ref{Prop1}.
Actually we shall prove more precisely that
$$  \vc{\widetilde N}_{\ep,\delta} \:  =  \:
F_{\eps, \delta} \, \na_h \phi \: + \: g_{\ep,\delta}  \na_h
\tilde \vr  \: + \:  \vc{s}_{\ep,\delta}
$$
for   explicit  functions~$\phi= \phi(t,x_h)$, $F_{\ep, \delta} = F_{\eps, \delta}(\tilde \vr)$,
 $ g_{\ep,\delta} = g_{\ep,\delta}(t, x_h)$  and a remainder $\vc s_{\ep,\delta} = \vc s_{\ep,\delta}(t,x_h)$
 satisfying
 {
 \begin{equation} \label{sepdelta}
 \lim_{\delta \rightarrow 0} \left( \limsup_{\eps \rightarrow 0} \: | \left( \vc s_{\ep,\delta} | \varphi \right) | \right) = 0
 \: \mbox{ for all } \: \varphi\in C^\infty_c ( (0,T)\times\Omega; R^2).
 \end{equation}
 }
 This  will of course  imply the result.

 \medskip

  {\em The notation $o(1)$  will refer hereafter to any term
  $\: s_{\ep,\delta} \: $satisfying \eqref{sepdelta}. Moreover     we  drop from now on the lowerscript~$\delta$, except if some ambiguity is liable to occur. }

 \medskip

Let us set
 $$ \vc{\widetilde N}_\eps  = \vc{\widetilde N}_\eps^1 + \vc{\widetilde N}_\eps^2 \: \equiv \:
 \Divh \left( \frac{1}{\tilde \vr} \va{\vc V_{\eps ,h}} \otimes \va{\vc V_{\eps,h}} \right) \: + \:
 \Divh \left(\frac{1}{\tilde \vr} \va{\{ \vc V_{\eps,h} \} \otimes \{ \vc V_{\eps,h} \} } \right), $$
{ where the notation $\va{\cdot}$, $\Bl \cdot \Br$ has been
introduced in
 (\ref{[v]}), and  treat these  parts separately.}
 \begin{itemize}
 \item {\bf Treatment of $\vc{\widetilde N}_\eps^1$:}
 \end{itemize}
 %We write $V = \tilde \vr \,  u$.
We notice that
 \begin{align*}
 \vc{\widetilde N}_\eps^1 & \: = \:   \Divh \left(    \, \va{\vc V_{\eps ,h}} \otimes \frac{1}{ \tilde \vr }
 \va{\vc V_{\eps ,h}} \right) \: = \:  \frac{1}{ \tilde \vr } \va{\vc V_{\eps ,h}} \Divh(\, \va{\vc V_{\eps ,h}})
 \: + \:  \va{\vc V_{\eps ,h}} \cdot \na_h  \left( \frac{1}{ \tilde \vr }   \va{\vc V_{\eps ,h}}  \right)\\
 & \: =  \:
  \frac{1}{ \tilde \vr }\va{\vc V_{\eps ,h}}  \Divh(\, \va{\vc V_{\eps ,h}})
  \: + \:   \frac{1}{2} \tilde \vr \, \na_h \left( | \frac1{\tilde \vr}\va{\vc V_{\eps ,h}} |^2\right) \:
  + \:      \big({\rm curl}_h \frac1{\tilde \vr}\va{\vc V_{\eps ,h}} \big) \va{\vc V_{\eps ,h}}^\perp   .
 \end{align*}
 On the one hand, averaging  \eqref{NSbisreg1} with respect to $x_3$, and multiplying by $\va{\vc V_{\eps ,h}}$,
 we get that
 $$  \frac{1}{ \tilde \vr }\va{\vc V_{\eps ,h}}  \Divh( \va{\vc V_{\eps ,h}}) = - \eps (\pa_t \va{r_\ep})
 \frac{1}{ \tilde \vr }\va{\vc V_{\eps ,h}} = - \eps \pa_t \left( \frac{1}{ \tilde \vr }\va{r_\ep}
 \va{\vc V_{\eps ,h}} \right) \: + \: \frac1 {\tilde \vr} \va {r_\ep} \, \eps \left(\pa_t  \va{\vc V_{\eps ,h}}\right), $$
where the first term at the right-hand side is of order $o(1)$. As regards the second term,
averaging the horizontal components of \eqref{NSbisreg2} with respect to $x_3$ and multiplying by $\va{r_\ep}$,
we end up with
\begin{align*}
  \frac1{\tilde \vr}\va{r_\ep} \, \eps \left(\pa_t  \va{\vc V_{\eps ,h}}\right)&  =
  - \na_h \left( P'(\tilde \vr) \,  \va{r_\ep} \right)   \va{r_\ep} \: - \:
  \frac1{\tilde \vr}\va{\vc V_{\eps ,h}}^\perp \va{r_\ep} \: + \:   \frac1{\tilde \vr}
\va{\eps   \vc{F}^{1}_{\ep,h}  +   \vc{F}^{2}_{\ep,h} }\va{r_\ep} \\
&  = \: - \frac{  1}{P'(\tilde \vr)} \na_h \left( \frac{1}{2}
\left|P'(\tilde \vr) \va{r_\ep} \right|^2 \right)  \: - \:
\frac1{\tilde \vr}\va{\vc V_{\eps ,h}}^\perp   \va{r_\ep}   \: +
\: o(1),
\end{align*}
where the $o(1)$ comes from the properties of the
$  \vc{F}_\ep^i$'s and the fact that $r_\ep = r_{\eps,\delta}$ is
uniformly bounded in $L^2_{loc} ((0,T)\times\Omega) $ with respect
to $\ep$ and $\delta$, see \eqref{weakconvrho}.
 Thus  $\vc{\widetilde N}_\eps^1$ can be written in the form
$$ \vc{\widetilde N}_\eps^1 =    \left(  \Bigl({\rm curl}_h \frac1{\tilde \vr}\va{\vc V_{\eps ,h}} \Bigr)
-    \frac1{\tilde \vr} \va{r_\ep}  \right)\va{\vc V_{\eps
,h}}^\perp \: + \: F(\tilde \vr) \, \na_h \phi  \: + \: o(1). $$
It remains to handle the first term at the r.h.s. {Therefore,
we average \eqref{NSbisreg2} with respect to $x_3$, divide by
$\tilde\vr$, take the curl of the horizontal components, and
subtract average of \eqref{NSbisreg1} with respect to $x_3$. We
obtain,
\begin{equation} \label{equation2Dr}
 \eps \pa_t  \left(    \Bigl({\rm curl}_h  \frac1{\tilde \vr}\va{\vc V_{\eps ,h}}  \Bigr)
 -  \frac1{\tilde \vr}\va{r_\ep} \right) \: = \:  \frac1{\tilde \vr^2}\va{\vc V_{\eps ,h}}
 \cdot \na_h \tilde \vr + \widetilde F_\ep^h
 \end{equation}
 }
 where
 $$
 \widetilde F_\ep^h :={\rm curl}_h  \frac1{\tilde \vr} \va{\ep \vc{F}_{\ep,h}^{1} + \vc{F}_{\ep,h}^{2}  }.
 $$
 According to Proposition~\ref{Prop2} there is a function~$f(\delta)$ going to zero with~$\delta$ such that, uniformly in~$\ep$,
  \begin{equation}\label{deffdelta}
  \sup_{\ep > 0}
 \left\|{\rm curl}_h  \frac1{\tilde \vr}  \vc{F}_{\ep,h}^{2}\right\|_{L^2((0,T) \times K;R^2)} \leq f(\delta).
  \end{equation}
 Now let us notice that by definition of~$ \tilde \vr$ one has
 $$
 \na_h \tilde \vr  \:  \: P'( \tilde \vr) = 2x_h, \quad \mbox{so} \quad
 \na_h  \tilde \vr \:  \:  p'(   \tilde \vr) = 2 \tilde \vr \: x_h.
 $$
 By Assumption~(\ref{HYP1}),  this
 gives in particular that there is a constant~$C$  such that for all~$x_h$,
  \begin{equation}\label{lowerboundtildevr}
 | \na_h  \tilde \vr ( x_h)| \geq C | x_h|.
 \end{equation}
Now let us consider a smooth function~$\chi_\delta$ defined by
$$\chi_\delta (x_h ) := \chi \left(\frac{\na_h \tilde \vr(x_h )}{ \sqrt{f(\delta)}} \right) \cdotp
$$
where  $\chi$ is  a function
of $ \DC ( R^2 ;[0,1])$ such that $\chi(x_h)=1$ if $|x_h|\leq 1$.
Using~(\ref{estimateV1}) we get
% \begin{align*}
%&\chi_\delta  \left(    \Bigl({\rm curl}_h  \frac1{\tilde \vr}\va{V_{\eps ,h}}  \Bigr)    -  \frac1{\tilde \vr}\va{r_\ep} \right)\va{V_{\eps ,h}}^\perp  \\
%& =   \:  \chi_\delta  \left(  \eps T_{\ep,\delta} ^1 + T_{\ep,\delta} ^2  -  \frac1{\tilde \vr}\va{r_\ep} \right) \left( \eps t_{\ep,\delta} ^1 + t_{\ep,\delta} ^2 \right) .
%\end{align*}
{
$$
\chi_\delta  \left(    \Bigl({\rm curl}_h  \frac1{\tilde
\vr}\va{\vc V_{\eps ,h}}  \Bigr) -  \frac1{\tilde \vr}\va{r_\ep}
\right)\va{\vc V_{\eps ,h}}^\perp   =
 \:  \chi_\delta  \left(  \eps T_{\ep,3} ^1 + T_{\ep,3} ^2  -
 \frac1{\tilde \vr}\va{r_\ep} \right) \va{ \eps \vc t_{\ep,h} ^1 + \vc t_{\ep,h} ^2 } ^\perp.
$$
}
Now we can write, by H\"older's inequality and the continuous
embedding~$W^{1,2} \subset L^6 $
 \begin{align*}
&\left\|\chi_\delta  \left(  \eps T_{\ep,3} ^1 + T_{\ep,3} ^2  -  \frac1{\tilde \vr}\va{r_\ep} \right)
\va{ \eps \vc t_{\ep,h} ^1 + \vc t_{\ep,h} ^2}^\perp  \right\|_{L^1((0,T) \times K;R^2 )} \\
& \leq \|\chi_\delta \|_{L^3(K)} \left\| \eps T_{\ep,3} ^1 +
T_{\ep,3} ^2  -  \frac1{\tilde \vr}\va{r_\ep}
\right\|_{L^{2}(0,T;L^2(K))} \| \eps \vc t_{\ep,h} ^1 + \vc
t_{\ep,h} ^2\|_{L^2(0,T;W^{1,2}(K;R^2))}
\end{align*}
Recalling that~$r_\ep$ is uniformly bounded  in~$\ep,\delta$
in~$L^\infty(0,T;L^2_{\rm loc}(\Omega))$, and using~(\ref{estimatecurlV2})
and~(\ref{lowerboundtildevr})  we infer that
 \begin{align*}
\left\|\chi_\delta  \left(  \eps T_{\ep,3} ^1 + T_{\ep,3} ^2  -  \frac1{\tilde \vr}\va{r_\ep} \right)
\va{ \eps \vc t_{\ep,h} ^1 + \vc t_{\ep,h} ^2}^\perp  \right\|_{L^1((0,T) \times K;R^2)} &\leq (\ep c(\delta) + c) \left|
\Big\{x_h, \:  |\na_h \tilde \vr(x_h ) |\leq \sqrt{f(\delta)}\Big\}  \right|^\frac13 \\
&\leq (\ep c(\delta) + c) \left| B(0,C  \sqrt{f(\delta)})  \right|^\frac13 \\
& = o(1).
\end{align*}
This allows to conclude that
$$
\chi_\delta  \left(    \Bigl({\rm curl}_h  \frac1{\tilde
\vr}\va{\vc V_{\eps ,h}}  \Bigr)    -  \frac1{\tilde
\vr}\va{r_\ep} \right)\va{\vc V_{\eps ,h}}^\perp = o(1).
$$

Then we can write
\begin{align*}
& (1- \chi_\delta  ) \left(    \Bigl({\rm curl}_h  \frac1{\tilde \vr}\va{\vc V_{\eps ,h}}  \Bigr)-
\frac1{\tilde \vr}\va{r_\ep} \right)\va{\vc V_{\eps ,h}}^\perp  \\
& = \:  (1- \chi_\delta  )  \left(    \Bigl({\rm curl}_h  \frac1{\tilde \vr}\va{\vc V_{\eps ,h}}  \Bigr)-
\frac1{\tilde \vr}\va{r_\ep} \right) \left(  \frac{\va{\vc V_{\eps ,h}}^\perp \cdot \na_h^\perp \tilde \vr}
{|\na_h \tilde \vr|^2}   \: \na_h^\perp \tilde \vr \: + \:   \frac{\va{\vc V_{\eps ,h}}^\perp \cdot
\na_h  \tilde \vr}{|\na_h \tilde \vr|^2} \:  \na_h  \tilde \vr \right) \\
& =  (1- \chi_\delta  )    \left(    \Bigl({\rm curl}_h
\frac1{\tilde \vr}\va{\vc V_{\eps ,h}}  \Bigr)-  \frac1{\tilde
\vr}\va{r_\ep} \right)  \frac{\va{\vc V_{\eps ,h}} \cdot \na_h
\tilde \vr}{|\na_h \tilde \vr|^2} \:  \na_h^\perp \tilde \vr  \:\:
+ \: g^1 \,   \na_h  \tilde \vr.
\end{align*}
for
$$g^1 \: \equiv \:    (1- \chi_\delta  )
\left(    \Bigl({\rm curl}_h  \frac1{\tilde \vr}\va{\vc V_{\eps
,h}} \Bigr)-  \frac1{\tilde \vr}\va{r_\ep} \right)  \,
\frac{\va{\vc V_{\eps ,h}}^\perp \cdot \na_h  \tilde \vr}{|\na_h
\tilde \vr|^2} \cdotp $$
 As regards the first term at the right-hand side of the last equality, we use \eqref{equation2Dr} to derive
\begin{align*}
&  { (1 - \chi_\delta)}  \left(    \Bigl({\rm curl}_h  \frac1{\tilde \vr}\va{\vc V_{\eps ,h}}  \Bigr)-
\frac1{\tilde \vr}\va{r_\ep} \right) \frac{\va{\vc V_{\eps ,h}} \cdot \na_h \tilde \vr}{|\na_h \tilde \vr|^2} \:
\na_h^\perp \tilde \vr \\
  =&  -\eps \pa_t \left(  \frac{   (1- \chi_\delta  ) }{2 |\na_h \tilde \vr|^2} \left( \tilde \vr
  \,\Bigl({\rm curl}_h \frac1{\tilde \vr}\va{\vc V_{\eps ,h}} \Bigr) -   \va{r_\ep}  \right) ^2 \na_h^\perp
  \tilde \vr \right)   +
  \tilde \vr^2  \left(  \eps T_{\ep,3} ^1 + T_{\ep,3} ^2  -  \frac1{\tilde \vr}\va{r_\ep} \right)
  \widetilde F_\ep^h \frac{\na_h^\perp \tilde \vr}{|\na_h \tilde \vr|^2}  (1- \chi_\delta  )
   \end{align*}
 with the notation introduced in~(\ref{estimateV1}).
We get, by~(\ref{estimatecurlV2}),
\begin{align*}
\left\|  \left(   \eps T_{\ep,3} ^1 + T_{\ep,3} ^2  -
\frac1{\tilde \vr}\va{r_\ep} \right) \widetilde F_\ep^h
\frac{\na_h^\perp \tilde \vr}{|\na_h \tilde \vr|^2}  (1-
\chi_\delta  )\right\|_{L^2(0,T; L^1(K))} & \leq  \left\| \frac{
(1- \chi_\delta  )}{|\na_h \tilde \vr|} \right\|_{ L^\infty(K)}
(\ep c(\delta)  + c) \: \|\widetilde F_\ep^h\|_{L^2((0,T) \times
K)}
\\
 & \leq \ep c(\delta)  + \frac{c}{\sqrt{f(\delta)}}\| \widetilde F_\ep^h \|_{L^2((0,T) \times K)}
\end{align*}
which is~$o(1)$ thanks to~(\ref{deffdelta}).
 Combining the previous inequalities leads to
\begin{equation*}
\vc{\widetilde N}_\eps^1=  g^1 \,   \na_h  \tilde \vr \: + \: F(\tilde \vr) \, \na_h \phi  \: + \: o(1),
 \end{equation*}
as expected.

\begin{itemize}
 \item {\bf Treatment of $\vc{\widetilde N}_\ep^2$:}
 \end{itemize}

   This time, we consider
 \begin{align*}
 \vc{\widetilde N}_\ep^2 & \: = \:
 \Divh \left(\frac{1}{\tilde \vr} \va{\Bl\vc V_{\eps ,h}\Br \otimes \Bl\vc V_{\eps ,h}\Br } \right) \\
 & \: = \: \left\langle \frac{1}{\tilde \vr} \,  \Bl\vc V_{\eps ,h}\Br \Divh \left(\Bl\vc{V}_{\eps ,h}\Br\right)
  \right\rangle \: + \:  \frac{1}{2} \tilde \vr \,  \left\langle \na_h \left| \frac{1}{\tilde \vr}
  \Bl\vc V_{\eps ,h}\Br \right|^2  \right\rangle\: + \:  \,
  \left\langle {\rm curl}_h \Bigl( \frac{1}{\tilde \vr} \Bl\vc V_{\eps ,h}\Br \Bigr)\,  \Bl\vc V_{\eps ,h}\Br^\perp \right\rangle.
 \end{align*}
As $\displaystyle  \frac{1}{\tilde \vr} \Bl\vc V_{\eps ,h}\Br$ has
zero vertical average, we can write
 $$   \curl \frac{1}{\tilde \vr} \Bl\vc V_{\eps }\Br \: = \: \left( \begin{matrix} \pa_3 \vc\Omega_{\ep,h} \\
 \omega_\ep \end{matrix} \right), \quad
 \vc\Omega_{\ep,h} \: :=  \: \frac{1}{\tilde \vr} \Bl\vc V_{\eps ,h}\Br^\perp -\pa_3^{-1} \na_h^{\perp}
 \frac{1}{\tilde \vr} \Bl\vc V_{\eps ,3}\Br, \quad \omega_\ep \: := \:  {\rm curl}_h  \frac{1}{\tilde \vr} \Bl\vc V_{\eps ,h}\Br,
 $$
{where we have set $\partial^{-1}_{3} a = I(a)$, see
(\ref{[v]}).} Applying ${\bf curl}_x$ to the momentum equation
\eqref{NSbisreg2} yields
\begin{align}
\label{eqcurl1}
\eps \pa_t \vc\Omega_{\ep,h}  \: & = \: \frac{1}{\tilde \vr} \Bl\vc V_{\eps ,h}\Br\: + \:   \Bl\pa^{-1}_3
({\bf curl}_x  \frac{1}{\tilde \vr}(\eps   \vc{F}_\eps^1 +   \vc{F}_\eps^2))_h \Br, \\
\label{eqcurl2} \eps \pa_t \omega_\ep \: & = \:  - \Divh
\frac{1}{\tilde \vr} \Bl \vc V_{\eps ,h}\Br\  \: + \: \left\{
{\rm curl}_h \frac{1}{\tilde \vr} (\eps   \vc{F}_\eps^1 +
\vc{F}_\eps^2)_h \right\}.
\end{align}
From there, we deduce that
\begin{align*}
   \left\langle {\rm curl}_h \Bigl( \frac{1}{\tilde \vr} \Bl\vc V_{\eps ,h}\Br \Bigr)\,
   \Bl\vc V_{\eps ,h}\Br^\perp \right\rangle & = \:   \left\langle   \omega_\ep \,   \Bl\vc V_{\eps ,h}\Br^\perp  \right\rangle
  = \: \tilde \vr \left\langle \omega_\eps \,  \eps \pa_t \vc \Omega_{\ep,h} ^\perp  \right\rangle \: + \: o(1)  \\
&  = \,  -\tilde \vr  \left\langle  (\eps \pa_t \omega_\ep)  \vc
\Omega_{\ep,h} ^\perp  \right\rangle \: + \:
 \left\langle \eps \pa_t \left( \tilde \vr \, \omega_\ep \,  \vc \Omega_{\ep,h} ^\perp \right)  \right\rangle  \: + \: o(1) \\
 & =    \, \tilde \vr \left\langle  \vc \Omega_{\ep,h} ^\perp  \Divh  \frac{1}{\tilde \vr} \Bl\vc V_{\eps ,h}\Br  \right\rangle  \: + \: o(1)
 \end{align*}
 Here, we have used repeatedly  equations (\ref{eqcurl1}),
 (\ref{eqcurl2}).
 %on $\Omega_{\ep,h}$ and $\omega_{\ep}$.
 The error terms generated by
 $$
  \tilde \vr \pa^{-1}_3 {\bf curl}_x  \frac{1}{\tilde \vr} \va{ \eps   \vc{F}_\ep^1 + \vc{\tilde {F}}_\ep^2}
 $$ are responsible for the $o(1)$ term, as can be verified using as previously~(\ref{estimateV1})-(\ref{estimatecurlV2}).
 Therefore
  \begin{align*}
 \vc{\widetilde N}_\ep^2 & \: = \:   \left\langle \frac{1}{\tilde \vr} \,  \Bl\vc V_{\eps ,h}\Br
 \Divh \left(\Bl\vc V_{\eps ,h}\Br\right)  \right\rangle \: +\:  \tilde \vr
 \left\langle \vc \Omega_{\ep,h} ^\perp
 \Divh  \frac{1}{\tilde \vr} \Bl\vc V_{\eps ,h}\Br  \right\rangle    \: + \:  F(\tilde \vr)  \na_h \psi    \: + \: o(1)  \\
 & \: = \:  \frac1{\tilde \vr } \, \left\langle\left( \Bl\vc V_{\eps ,h}\Br
 + \tilde \vr\vc \Omega_{\ep,h} ^\perp \right)\,  \Divh  \left(\Bl\vc V_{\eps ,h}\Br\right)
 \right\rangle  \: + \:  \left\langle  \tilde \vr \vc \Omega_{\ep,h} ^\perp   \Bl\vc V_{\eps ,h}\Br
 \cdot \na_h \frac1 {\tilde \vr}   \right\rangle   \: + \:  F(\tilde \vr)  \na_h \psi    \: + \: o(1) .
 \end{align*}
% Whence
%
 By straightforward manipulations we have
 \begin{equation*}
 \frac1 {\tilde \vr}\left( \Bl\vc V_{\eps ,h}\Br  + \tilde \vr\vc \Omega_{\ep,h} ^\perp \right)  \pa_3 \Bl V_{\eps ,3}\Br
 \:   = \:  \pa_3  \left(
  \frac1 {\tilde \vr}  \Bl V_{\eps ,3}\Br  \left( \Bl\vc V_{\eps ,h}\Br  + \tilde \vr\vc \Omega_{\ep,h} ^\perp \right) \right)
   \: - \:  \frac{1}{2 \tilde \vr} \na_h \big | \Bl V_{\eps ,3} \Br \big |^2  -  \big | \Bl V_{\eps ,3} \Br \big |^2
   \na_h \frac1 {\tilde \vr}
 \end{equation*}
 %so that $\widetilde N_\ep^2$ has the following equivalent form, for some function~$g^2$:
 %$$ \widetilde N_\ep^2\: = \:   \frac1 {\tilde \vr}    \left\langle \left( \Bl\vc V_{\eps ,h}\Br\:
 %+  \tilde \vr\vc\Omega_{\ep,h} ^\perp \right) \Div \Bl\vc V_{\eps }\Br
 % \right\rangle  \:
 %+ \:  \left\langle  \tilde \vr\Omega_{\ep,h} ^\perp   \BlV_{\eps ,h}\Br \cdot \na_h \frac1 {\tilde \vr}\right\rangle    \:
 %+ \:  F(\tilde \vr)  \na_h \psi    \: + \: g^2    \na_h  {\tilde \vr}  \: + \: o(1).   $$

 We can replace the first term by using equation \eqref{NSbisreg2}
 $$ \frac1 {\tilde \vr}    \left\langle  \left( \Bl\vc V_{\eps ,h}\Br +  \tilde \vr \vc \Omega_{\ep,h} ^\perp \right)
 \Div \Bl\vc V_{\eps }\Br  \right\rangle =   - \frac1 {\tilde \vr}  \left\langle  ( \eps \pa_t \Bl r_\ep\Br   )
 \left( \Bl\vc V_{\eps ,h}\Br +  \tilde \vr\vc\Omega_{\ep,h} ^\perp \right) \right\rangle $$
 which leaves us with
 \begin{align*}
 \vc{\widetilde N}_\ep^2 \: & = \:  - \frac1 {\tilde \vr}  \left\langle  ( \eps \pa_t \Bl r_\ep\Br   )
  \left( \Bl\vc V_{\eps ,h}\Br +  \tilde \vr\vc \Omega_{\ep,h} ^\perp \right) \right\rangle\: + \:
  \left\langle  \tilde \vr\vc\Omega_{\ep,h} ^\perp   \Bl\vc V_{\eps ,h}\Br \cdot \na_h \frac1 {\tilde \vr}
  \right\rangle    \: + \:  F(\tilde \vr)  \na_h \psi    \: + \: g^2    \na_h  {\tilde \vr}  \: + \: o(1)    \\
& =  \frac1 {\tilde \vr}  \left\langle \Bl r_\ep \Br   \,  \eps \pa_t
\left( \Bl\vc V_{\eps ,h}\Br +  \tilde \vr\vc\Omega_{\ep,h} ^\perp
\right) \right\rangle \: + \:  \left\langle  \tilde
\vr\vc\Omega_{\ep,h} ^\perp   \Bl\vc V_{\eps ,h}\Br \cdot \na_h \frac1
{\tilde \vr} \right\rangle    \: + \:  F(\tilde \vr)  \na_h \psi
\: + \: g^2 \na_h  {\tilde \vr}  \: + \: o(1) .
\end{align*}
Now, thanks to \eqref{NSbisreg2} and \eqref{eqcurl1},
$$ \frac1 {\tilde \vr}  \left\langle \Bl r_\ep\Br  \,  \eps \pa_t \left( \Bl\vc V_{\eps ,h}\Br +  \tilde \vr
\vc \Omega_{\ep,h} ^\perp \right) \right\rangle \: = \:     -\frac{1}{2 P'(\tilde \vr)} \na_h \left| P'(\tilde \vr) \Bl r_\ep \Br \right|^2 + o(1), $$
so that $\vc{\widetilde N}_\ep^2 $ resumes to
$$ \vc{\widetilde N}_\ep^2 \: =  \:   \left\langle  \tilde \vr\vc \Omega_{\ep,h} ^\perp   [\vc V_{\eps ,h}]
\cdot \na_h \frac1 {\tilde \vr}   \right\rangle    \: + \:  F(\tilde \vr)  \na_h \psi    \: + \: g^2    \na_h  {\tilde \vr}  \: + \: o(1) .$$
Finally, we proceed as previously in the case of $\vc{\widetilde N}_\ep^1$, with the first term at the right-hand side, this time omitting  for simplicity the cut-off near~$\na_h \tilde \vr = 0$. We write
\begin{align*}
 \left( \Bl \vc V_{\eps ,h} \Br \cdot \na_h \tilde \vr  \right)\vc \Omega_{\ep,h} ^\perp
&  = \:  \left( \Bl \vc V_{\eps ,h} \Br \cdot \na_h \tilde \vr  \right) \left( \frac{\vc \Omega_{\ep,h}   \cdot
\na_h \tilde \vr}{|\na_h \tilde \vr|^2}  \na_h^\perp \tilde \vr \: + \:  \frac{\vc \Omega_{\ep,h} ^\perp \cdot
\na_h  \tilde \vr}{|\na_h \tilde \vr|^2}    \na_h  \tilde \vr \right) \\
& = \tilde \vr \left(\eps \pa_t\vc \Omega_{\ep,h}  \cdot \na_h
\tilde \vr  \right) \frac{\vc \Omega_{\ep,h}   \cdot \na_h \tilde
\vr}{|\na_h \tilde \vr|^2}
 \na_h^\perp \tilde \vr  \: + \: g^3 \na_h  \tilde \vr \: + \: o(1)
\end{align*}
where we have used \eqref{eqcurl1} { in the passage} from the
second to the third line. This yields
\begin{align*}
\left( \Bl \vc V_{\eps ,h} \Br \cdot \na_h \tilde \vr  \right)
\vc\Omega_{\ep,h} ^\perp  =  g \na  \tilde \vr \: + \: o(1)
\end{align*}
and finally we get an expression of the form:
$$ \vc{\widetilde N}_\ep^2 \: =¬†\: F(\tilde \vr)\na_h \psi \: + \: g \na_h \tilde \vr  \: + \: o(1) $$
as expected. Proposition \ref{Prop1} is proved. 

\qed

\subsection{Regularization process} \label{subsecregular}
 In this section we shall prove Proposition~\ref{Prop2}.
 We start by recalling that
 $$
  \eps \pa_t r_{\eps} + \Div (\vc V_{\ep})   = 0,
  $$
   and
$$
  \eps \pa_t \vc V_{\ep} +  \tilde \vr \Grad \left( P'(\tilde \vr) \,  r_{\eps}  \right) + \vc{b} \times \vc V_{\ep}  =
\eps    \vc{N}_{\eps}+  \eps    \vc{F}_{\eps} ,
$$
with
 notation~(\ref{defNep}) and~(\ref{defSep}).
 The first step consists in establishing some bounds for~$  \vc{N}_\ep$ and $   \vc{F}_\ep$. The energy bound clearly
 implies that~$\vc{N}_\ep$ is bounded in~$L^\infty(0,T;W^{-1,1} (K;R^3)) $. As for~$\vc{F}_\ep$, one has clearly that~$
\Div \tn{S}(\Grad \vue) $ is bounded in~$L^2(0,T;W^{-1,2}
(K;R^3)) $, and~$\displaystyle \frac{1}{\eps^2} \Grad \left(
p(\vre) - p(\tilde \vr) -  p'(\tilde \vr) (\vre - \tilde
\vr)\right) $ is bounded in~$L^\infty(0,T;W^{-1,1} (K;R^3)) $.
  Therefore in particular
  \begin{equation}\label{bddnepsep}
  \vc{N}_\ep + \vc{S}_\ep \quad \mbox{ is bounded in†}† \quad L^2(0,T; W^{-5/2,2}(K;R^3))
  \end{equation}
for any compact $K \subset \Omega$.

\medskip
\noindent
Now let us proceed to the regularization.  First we notice that
$$
  {\vc V_{\eps , \delta}} = \eps  \kappa_\delta * (r_\eps \, \vue)    \: +
  \:   \kappa_\delta * (\tilde \vr \vue)   \:  = :\: \eps \vc t_{\ep,\delta} ^1 +\vc t_{\ep,\delta} ^2
$$
and
$$
  {\bf curl}_x \left( \frac1{\tilde \vr}\ \vc V_{\eps , \delta} \right )=
   \eps  {\bf curl}_x \left( \frac{1}{\tilde \vr} \kappa_\delta * (r_\eps \, \vue)  \right)   \: + \:
    {\bf curl}_x \left( \frac{1}{\tilde \vr} \kappa_\delta * (\tilde \vr \vue)  \right) \:  =: \: \eps \vc T_{\ep,\delta} ^1 +
    \vc T_{\ep,\delta} ^2
$$
so thanks to the $L^2 $ bound on $r_\eps$ and the $W^{1,2}$ bound on $\vu_\eps$, we deduce easily that for all $k,K$
 $$
{\| \vc t_{\ep,\delta} ^1 \|_{L^2(0,T; W^{k,2}(K;R^3))} + \| \vc
T_{\ep,\delta} ^1 \|_{L^2(0,T; W^{k,2}(K;R^3))} }\:\le \: c(\delta),
$$
$$
\quad \| \vc t_{\ep,\delta} ^2 \|_{L^2(0,T ; W^{1,2}( K;R^3))} + \|
\vc T_{\ep,\delta} ^2 \|_{L^2((0,T) \times K;R^3)} \le c  \:
$$
uniformly in $\eps$ (and $\delta$ for the second bound).  This proves~(\ref{estimatecurlV2}).
The uniform bounds derived previously also
give directly the
convergences (\ref{delta-cv}).

\medskip
\noindent
Now let us turn to the wave equations.
By convolution we get (with obvious notation)
\bFormula{nader1}
 \ep \partial_t r_{\ep,\delta} +  \Div\vc  V_{\ep,\delta} = 0,
\eF
and
\bFormula{nader2}
\eps \pa_t \vc V_{\ep,\delta} + \ \tilde \vr \, \Grad \left(
P'(\tilde \vr) \,  r_{\eps,\delta}  \right) + \vc{b} \times \vc
V_{\ep,\delta}  = \eps  \vc{F}^1_{\eps, \delta}  +  \vc{F}^2_{\eps, \delta}
\eF
with
$$
 \vc{F}^1_{\ep, \delta} = \vc{N}_{\ep, \delta}  + \vc{F}_{\ep, \delta}
$$
and
$$
  \vc{F}^2_{\eps, \delta} =    \tilde \vr \, \Grad \left( P'(\tilde \vr) \,  r_{\eps,\delta}  \right) - \left( \tilde \vr \, \Grad \left( P'(\tilde \vr) \,  r_{\eps}  \right) \right) * \kappa_\delta .
$$
 Clearly~(\ref{bddnepsep}) implies~(\ref{bddS1epsdelta}), so let us turn to the statement (\ref{bddS2epsdelta}).

\medskip
\noindent
In order to see (\ref{bddS2epsdelta}),
we use \cite[Proposition 4.1]{masmoudi} (which forces the restriction~$\gamma >3$) on compactness
of solutions to (\ref{nader1}), (\ref{nader2}), namely,
\[
\| r_\eps-r_{\ep,\delta}\|_{L^p([0,T]; L^2(K))} \to 0 \hbox{ as } \delta \to 0
\ \mbox{for any compact}\ K \subset \Omega \ \mbox{and any} \ p \geq 1,
\]
together with Lemma 3.3 (2) of~\cite{masmoudi}. Note that, compared with the situation
treated in \cite[Proposition 4.1]{masmoudi}, the present system contains an extra term
$\vc{b} \times \vc{V}_{\ep,\delta}$ already known to be compact with respect to the space variable.

\medskip
\noindent
The vanishing  of  $\vc{F}^2_{\ep,\delta}$ (uniformly in $\eps$)  follows directly.  To handle the convergence of  $\curl \frac{1}{\tilde \vr} \vc{F}^2_{\ep,\delta}$, we then  notice that
 $$ {\bf curl}_x \vc{F}^2_{\ep,\delta} = A(\tilde \vr) \na r_{\eps,\delta} - (A(\tilde \vr) \na r_\eps) * \kappa_\delta $$
for some smooth matrix function $A$. Still using Lemma 3.3 (2) in~\cite{masmoudi}, we obtain the vanishing of ${\bf curl}_x \vc{F}^2_{\ep,\delta}$, which together  with the one of  $\vc{F}^2_{\ep,\delta}$ completes the proof of Proposition \ref{Prop2}.

\qed

 \subsection{Conclusion}
Thanks to Proposition \ref{Prop2}, we can conclude the proof of
Theorem \ref{Tm2}. { We
keep the notation $\vc{N}_{\ep,h}$,~{$\vc{\widetilde N}_\eps$} and~$\vc{\widetilde N}_{\eps,\delta}$ of Section \ref{subsecconvective},
see (\ref{N}), (\ref{Nh}), (\ref{Nhd}).} 

\noindent 
Let~$\psi = \psi(t,x_h)
\in \DC((0,T) \times R^2)$ such that~$\na_h \tilde \vr \cdot
\na_h^\perp \psi = 0$. We write
\begin{align*}
&\left|  \left( \va{- \vc{N}_{\eps,h}} | \na_h^\perp \psi \right) - \left( \vc{\widetilde N}_{\eps,\delta} | \na_h^\perp
\psi \right)
\right|=\left|\int_0^T\int_\Omega\left(\vre \vu_{\ep,h}\otimes
\vu_{\ep,h}-\frac 1{\tilde\vr} \vc V_{\ep,\delta,h}\otimes \vc
V_{\ep,\delta,h}\right):\nabla_h \otimes \nabla^\perp_h \psi {\rm d}x{\rm
d}t\right|\\ \\
&  \:\le \:   \left|  \int_0^T \int_{\Omega}  \left( (\vc
V_{\eps,h} - \vc V_{\eps,\delta, h}) \otimes \vu_{\eps,h} \right)
 :  (\na_h \otimes \na_h^\perp \psi)\ {\rm d}x{\rm
d}t   \right| \\ \\
&
+  \left|  \int_0^T \int_{\Omega} \left( \vc V_{\eps,\delta,h} \otimes (\vu_{\eps,h} - \frac{\vc V_{\eps,h}}
 {\tilde \vr})
 \right) :  (\na_h \otimes \na_h^\perp \psi)\ {\rm d}x{\rm
d}t   \right|
\\ \\
 &  \:+ \: \left|  \int_0^T \int_{\Omega}    \left( \vc V_{\eps,\delta,h} \otimes  \frac{\vc V_{\eps,h} -
 \vc V_{\ep,\delta,h}}
 {\tilde \vr} \right) :  (\na_h \otimes \na_h^\perp \psi){\rm d}x{\rm d}t  \right|
  =: \: I^1_{\eps,\delta} + I^2_{\eps,\delta} + I^3_{\eps,\delta} \: .
\end{align*}
We have:
$$
 I^1_{\eps,\delta} =  \left|  \int_0^T \int_{\Omega}  \vc V_{\eps,h}  \cdot
 \left(  (\na_h \otimes \na_h^\perp \psi)\vu_{\eps,h} -
 \kappa_\delta * \big((\na_h \otimes \na_h^\perp \psi)\vu_{\eps,h}\big)  \right)  dx dt \right| = O(\delta)$$
 uniformly in $\eps$, using that
$$ \|   (\na_h \otimes \na_h^\perp \psi)\vu_{\eps,h} -    \kappa_\delta *  \big((\na_h \otimes
\na_h^\perp \psi)\vu_{\eps,h}\big) \|_{L^2((0,T) \times
\Omega;R^2)} \: \le \: C \delta \| \vu_{\eps,h} \|_{L^2(0,T;
W^{1,2}(K; R^2))}
$$ for some compact $K$ containing the support of $\psi$. 
Then, noticing that
$$ \frac{\vc V_\eps}{\tilde \vr} = \vu_\eps \: + \: \eps  \frac{r_\eps}{\tilde \vr} \vu_\eps $$
one obtains easily
$$ I^2_{\eps,\delta}  \le c (\delta)  \eps.  $$
Finally, we remark that
$$ \vc V_\eps - \vc V_{\ep,\delta} \: = \: \left(\kappa_\delta * (\tilde \vr \vu_\eps)- \tilde \vr \vu_\eps  \right)  \:
+ \:  \eps \left(  r_\eps \vu_\eps - \kappa_\delta * (r_\eps
\vu_\eps) \right) $$
to obtain
$$
 I^{3}_{\eps,\delta} \: \le \: c \, \delta  \: + \: c (\delta) \, \eps .
$$
Putting these inequalities altogether yields
$$ \limsup_{\delta \rightarrow 0} \, \limsup_{\eps \rightarrow 0} \: ( I^1_{\eps,\delta} + I^2_{\eps,\delta} + I^3_{\eps,\delta} ) \: = \: 0.
$$
Combining this with Proposition \ref{Prop2}, we deduce that
$$ \lim_{\eps \rightarrow 0} \left( \vc{N}_{\eps,h} | \na_h^\perp \psi \right)  \: = \:  0 $$
 for all $\psi = \psi(t,x_h) \in \DC ((0,T) \times R^2)$ such
that $\na_h \tilde \vr \cdot \na_h^\perp \psi =  0$, meaning the function $\psi$ is radially symmetric.

\medskip
\noindent
We are now at the point of getting the equation satisfied by  $r,\vu$, cf. (\ref{kernel1}).
The
horizontal part of the momentum equation reads
$$  \pa_t (\vre \vu_{\eps,h}) +  \Div \left( \vre \vu_{\ep} \otimes \vu_{\ep,h} \right)
+ \eps^{-2} \, \na_h p(\vre) + \eps^{-1}  \vre \vu_{\ep,h}^\perp =
[\Div \tn{S}(\Grad \vu_\ep)]_h  + \frac{\vre}{\eps^2} \, \na_h G .
$$
We recall that $\na_h G = \na_h P(\tilde \vr) = P'(\tilde \vr)
\na_h \tilde \vr$. We integrate with respect to $x_3$ the last
equation, and apply~${\rm curl}_h$. We obtain
\begin{equation} \label{curlhuh}
\begin{array}{r}
\pa_t {\rm curl}_h \va{\vre  \vu_{\eps,h}} \: + \:  {\rm curl}_h
\Divh \va{ \vre \vu_{\ep,h} \otimes \vu_{\ep,h} } + \eps^{-1}
\Divh\va{\vre \vu_{\ep,h}} \\ \\
=   \displaystyle   {\rm curl}_h \Divh \tn{S}_{h,h}(\nabla_h \va{\vu_{\ep,h}} )
+ {\rm curl}_h \left(\frac{ P'(\tilde \vr)
\langle\vre\rangle}{\ep^2}  \na_h \tilde \vr\right),
\end{array}
\end{equation}
where
$$
{\tn
S}_{h,h}(\nabla_h\vu_h)=\mu(\nabla_h\vu_h+\nabla_h^\perp\vu_h-\frac
23{\rm div}_h(\vu_h){\tn I}_{h,h})\;\mbox{with}\; \tn
I_{h,h}\;\mbox{identity matrix in}\; R^2.
$$
 Continuity equation yields $ \Divh\va{\vre  \vu_{\ep,h} } =
- \eps \pa_t \langle r_\eps \rangle$; we employ this fact and
\eqref{curlhuh}, where we use a radially symmetric test function $\psi \in \DC ((0,T) \times R^2)$ to get
\begin{align*}
&\Big(   \pa_t \big( {\rm curl}_h \va{\vre  \vu_{\eps,h}} -
\va{r_\eps}\big) -     {\rm curl}_h \Divh
\tn{S}_{h,h}(\Grad\langle\vu_{\ep,h}\rangle) |
\psi \Big) \\ \\
& \: = \: \left(  \langle N_{\ep,h}\rangle | \na^\perp_h \psi
\right) -  \frac{1}{\eps^2} \int_0^T
\int_{R^2} P'(\tilde \vr)
\langle\vre\rangle \na_h \tilde \vr \cdot \na_h^\perp \psi  \, dx_h  \, dt  \\
\end{align*}
 The first (convective) term at the r.h.s. goes to zero, whereas
the second one is identically zero by the properties of $\psi$.
All other quantities converge easily to yield
$$  \Big( \pa_t    \Big({\rm curl}_h \va{\tilde \vr \vc{U}_h} -   \va{ r }\Big) -
 \mu\Delta_h {\rm curl}_h \vc{U}_h,  \psi \Big) =
0  $$
for any radially symmetric $\psi$.  Using the properties  (\ref{kernel1}-- \ref{kernel2}), we arrive
at \eqref{limiteq}.

\medskip
Finally, repeating the same procedure with
$\psi  \in \DC ([0,T) \times R^2)$, $\psi$ radially symmetric, we obtain
\begin{equation} \label{limiteqbis}
   - \left(    ({\rm curl}_h (\tilde \vr \vu_h) -    r ) | \pa_t \psi
   \right) -   \mu\left(    \Delta_h    {\rm curl}_h  \vu_h |
   \psi \right) =   - \left(   {\rm curl}_h \langle \tilde \vr \vu_{0,h}\rangle -
   \langle r_0 \rangle |  \psi_{|t=0}  \right),
\end{equation}
where $(\vu_0, r_0)$ are the weak limits of the family of initial data
$\vu_{0,\eps}, r_{0,\eps}$. This justifies the initial condition stated in (\ref{INIT}). Moreover, in view of hypothesis (\ref{u1++}), we have
\begin{equation} \label{weightu0}
 \sqrt{P'(\tilde \vr)}  \langle r_0 \rangle \in L^2(R^2), \quad \sqrt{\tilde \vr} \langle \vu_{0,h} \rangle \in L^2(R^2; R^2).
 \end{equation}

Under these circumstances, it is easy to show that
\eqref{kernel1}-\eqref{kernel2}-\eqref{limiteqbis} admits a unique solution. Indeed, taking $\psi = P'(\tilde \vr) r(0)$ as a test function in \eqref{INIT}, we check that
\eqref{weightu0} implies
$$ \int_{R^2}\left(  \tilde \vr \left| \na_h \left( P'(\tilde \vr) r(0) \right) \right|^2 \: + \: P'(\tilde \vr) |r(0)|^2  \right) \, {\rm d}x \: < \: +\infty; $$
whence uniqueness follows from standard energy arguments. See   \cite[section 5.2]{FeGaNo}.
 Thus, there is exactly one accumulation point for the sequence $(r_\eps,
\vu_\eps)$, and the whole sequence converges to it.

\def\ocirc#1{\ifmmode\setbox0=\hbox{$#1$}\dimen0=\ht0 \advance\dimen0
  by1pt\rlap{\hbox to\wd0{\hss\raise\dimen0
  \hbox{\hskip.2em$\scriptscriptstyle\circ$}\hss}}#1\else {\accent"17 #1}\fi}

%\bibliography{citace}
%\bibliographystyle{plain}
\end{document}